\documentclass[onefignum,onetabnum,nohypdvips]{siamonline190516}

\usepackage{lipsum}
\usepackage{amsfonts}
\usepackage{graphicx}
\usepackage{epstopdf}
\usepackage{algorithmic}
\usepackage{enumerate}
\usepackage{csquotes}
\usepackage{todonotes}
\usepackage{subcaption}
\usepackage{wrapfig}
\usepackage{amssymb}
\usepackage{booktabs}
\ifpdf
  \DeclareGraphicsExtensions{.eps,.pdf,.png,.jpg}
\else
  \DeclareGraphicsExtensions{.eps}
\fi

\usepackage{soul}

\usepackage{enumitem}
\setlist[enumerate]{leftmargin=.5in}
\setlist[itemize]{leftmargin=.5in}

\newsiamremark{remark}{Remark}
\newsiamremark{hypothesis}{Hypothesis}
\crefname{hypothesis}{Hypothesis}{Hypotheses}
\newsiamthm{claim}{Claim}

\headers{An Inexact Semi-smooth Newton Method on Riemannian Manifolds}{W. Diepeveen, and J. Lellmann}

\title{An Inexact Semi-smooth Newton Method on Riemannian Manifolds with Application to Duality-based Total Variation Denoising}

\author{Willem Diepeveen\thanks{Department of Applied Mathematics and Theoretical Physics, University of Cambridge 
  (\email{wd292@cam.ac.uk}).}
\and Jan Lellmann\thanks{ Institute of Mathematics and Image Computing, University of Lübeck
  (\email{lellmann@mic.uni-luebeck.de}).}
}

\usepackage{amsopn}

\newcommand{\jla}[1]{#1}

\newcommand{\wdpa}[1]{#1}

\newcommand{\manopt}{\texttt{Manopt.jl}}

\ifpdf
\hypersetup{
  pdftitle={An Inexact Semi-smooth Newton Method on Riemannian Manifolds with Application to Duality-based Total Variation Denoising},
  pdfauthor={W. Diepeveen and J. Lellmann}
}
\fi

\begin{document}

\maketitle

\begin{abstract}
  We propose a \jla{higher-order} method for solving non-smooth optimization problems on manifolds. In order to obtain superlinear convergence, we apply a Riemannian Semi-smooth Newton method to a non-smooth non-linear primal-dual optimality system based on a recent extension of Fenchel duality theory to Riemannian manifolds. We also propose an inexact version of the Riemannian Semi-smooth Newton method and prove conditions for local linear and superlinear convergence that hold independent of the sign of the curvature. Numerical experiments on $\ell^2$-TV-like problems with dual regularization confirm superlinear convergence on manifolds with positive and negative curvature.
\end{abstract}



%

\section{Introduction}
Energy-based modeling in image- and data analysis requires the numerical minimization of large-scale energy functions. Due to the growing popularity of sparsity-based approaches such as compressive sensing \cite{donoho2006compressed} and total variation-based image processing \cite{rudin1992nonlinear}, these energies often incorporate non-smooth terms.

First-order methods based on (sub-)gradients for minimizing such energies have become very popular due to their robustness \cite{combettes2005signal,combettes2007douglas,yin2008bregman,wang2008new,chambolle2011first,bertsekas2011incremental}; see also the survey in \cite{esser2010general}. However, their convergence rate in the general case is typically limited to $\mathcal{O}(\frac{1}{\epsilon})$ iterations for achieving $\epsilon$-suboptimality, which can  sometimes be improved to $\mathcal{O}(\frac{1}{\sqrt{\epsilon}})$ using acceleration strategies \cite{he2014convergence,pock2011diagonal,wright2009sparse}.

Options for obtaining a superlinear convergence  rate include interior point methods \cite{fountoulakis2014matrix} or Newton-like methods, such as 
Quasi Newton \cite{yu2010quasi,chen2014fixed}, Proximal Newton \cite{becker2012quasi,pan2013sparse,lee2014proximal, byrd2016inexact,yue2019family}, Forward Backward Newton \cite{patrinos2014forward} and Semi-smooth Newton (SSN)  \cite{griesse2008semismooth, milzarek2014semismooth, byrd2016family, xiao2018regularized,li2018highly}. In this work we will focus on a generalization of the \wdpa{SSN} method, which was originally proposed in \cite{qi1993nonsmooth} and initially used in \cite{luca1996semismooth, sun1997newton}, but popularized through optimal control applications in the early 2000s \cite{Hintermueller2002}, before being discovered by the image processing community \cite{griesse2008semismooth,milzarek2014semismooth,byrd2016family}. We also refer to \cite{rust2017,xiao2018regularized} for some recent variants.

\begin{figure}[h!]
        \centering
\includegraphics[width=0.4\textwidth]{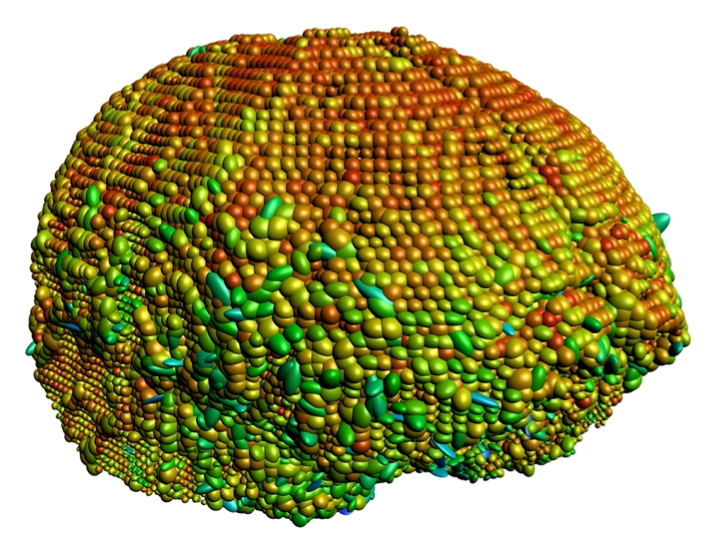}
\caption{Diffusion-Tensor Magnetic Resonance Imaging (DT-MRI) allows precise fiber tracing in the human brain by modeling the diffusion directions of water molecules using multivariate normal distributions, which are represented by elements on the manifold of symmetric positive definite matrices $\mathcal{P}(3)$. The resulting optimization problems are manifold-valued and often non-smooth, which necessitates specialized numerical optimization algorithms. \jla{Visualization courtesy of {\cite{Bergmann2016}}, based on the Camino data set {\cite{cook2006camino}}.}}\label{fig:camino}
\end{figure}

On the modeling side, there has been an increasing interest in manifold-valued data processing: apart from statistical \cite{pennec2006riemannian, pennec2006intrinsic, fletcher2007riemannian, laus2017nonlocal} and PDE approaches \cite{kimmel2002orientation, chefd2004regularizing} to smooth data processing on manifolds, non-smooth variational approaches on Riemannian manifolds have been gaining momentum as well in the 2010s. Typical applications include non-linear color spaces \cite{chan2001total} such as the Chromaticity Brightness model ($S^2 \times \mathbb{R}$) and the Hue Saturation Value model ($S^1 \times \mathbb{R}^2$), InSAR imaging \cite{massonnet1998radar} ($S^1$), \jla{Diffusion Tensor imaging {\cite{basser1994mr}} with the manifold of positive definite symmetric $3\times 3$ matrices $\mathcal{P}(3)$ (Fig.~{\ref{fig:camino}}) and Electron Backscatter Diffraction imaging {\cite{adams1993orientation}} with the 3D rotation group $SO(3)$.} Also several branches in data science have become more involved in using the geometry of problems. Examples include sparse principal component analysis, compressed mode analysis in physics, unsupervised feature selection and sparse blind convolution, where in these cases the Stiefel manifold captures all geometric information;  see \cite{chen2020proximal} for an overview.

In this work, we consider optimization problems of the form 
\begin{equation}
        \inf_{p \in \mathcal{M}} \left\{ F(p)+G(\Lambda(p)) \right\},
\end{equation}
where  $F: \mathcal{M} \rightarrow \overline{\mathbb{R}}$ and $G: \mathcal{N} \rightarrow \overline{\mathbb{R}}$ are non-smooth functions mapping into the extended real line $\overline{\mathbb{R}}:=\mathbb{R}\cup\{\pm\infty\}$, $\mathcal{M}$ and $\mathcal{N}$ are Riemannian manifolds, and  $\Lambda: \mathcal{M} \rightarrow \mathcal{N}$ is differentiable. \jla{In variational image processing models, $\mathcal{M}$ will typically be a discretized space of \emph{manifold-valued} functions, i.e., whose \emph{range} is restricted to a manifold. This should be contrasted with case where the \emph{domain} of the unknown function is a manifold, but the range is Euclidean or even scalar-valued, such as when computing quantities on the surface of a 2D or 3D shape. From an optimization viewpoint, the latter is somewhat easier, as -- after discretization -- the problem boils down to a finite-dimensional Euclidean optimization problem.}

Unlike the Euclidean case $\mathcal{M}=\mathbb{R}^{m}$, $\mathcal{N}=\mathbb{R}^{n}$, to our knowledge no intrinsic higher-order solvers with superlinear convergence exist in the manifold-valued case. Our goal is to close this gap.

\subsection{Related Work}
There have been multiple attempts to generalize algorithms for both smooth and non-smooth optimization on manifolds. Early work can be traced back to \cite{luenberger1972gradient}, where a manifold constraint is included in the optimization problem and enforced using projections. Using this so-called \emph{extrinsic} approach, algorithms from the real-valued case can be made to work on manifolds as well, if they can be reasonably embedded into en Euclidean space. While this technique does not capture the manifold structure specifically, 
recent contributions still tend to this approach for optimization \cite{lai2014splitting,bolte2014proximal,chen2016augmented,kovnatsky2016madmm,zhu2017nonconvex}. A step towards exploiting the  \emph{intrinsic} geometry of manifolds is to rely on local charts, which provide a linear subspace in which methods for real-valued optimization can be applied. However, both embedding and localization approaches suffer from serious drawbacks, which are discussed extensively in \cite[Sec. 1]{grohs2016nonsmooth}. Hence, intrinsic methods that do not rely on charts have become very popular.

Pioneering work using these ideas in smooth optimization on manifolds was done in the 1994 with \cite{smith1994optimization, udriste1994convex}, whose authors formulated several algorithms such as gradient descent, Newton's method and conjugate gradient on Riemannian manifolds. From then on, the community started working on generalizing other algorithms to Riemannian manifolds \cite{absil2007trust,bortoloti2018damped,castro2017secant}, specializing algorithms \cite{edelman1998geometry,abrudan2008steepest,he2013convergence} and the application to real-world problems \cite{adler2002newton,afsari2013convergence}. An extensive overview of first-and second-order methods for smooth optimization on matrix manifolds can be found in \cite{absil2009optimization}.

The development of non-smooth Riemannian optimization methods appears to have occurred somewhat independently in the image analysis- and in the non-smooth optimization communities. 
Whereas the former mostly relied on methods that generalize convexity, the latter has been focussed on finding new notions and restrictions to generalize manifold algorithms from the smooth to non-smooth setting. We give a short survey  of both, as we believe that either provide valuable insights.


\paragraph{Manifold-valued Imaging}
The non-smooth Rudin-Osher-Fatemi (ROF) model \cite{rudin1992nonlinear} is the prototype of modern image analysis. This model revolves around the notion of Total Variation (TV), which has also been generalized to the $S^1$ manifold in the early 1990s \cite{giaquinta1993variational}  and has been extended to general Riemannian manifolds in the second half of the 2000s \cite{giaquinta2006bv,giaquinta2007maps}. Although early numerical attempts for  TV-regularized problems on manifolds originate in the early 2000s \cite{chan2001total}, the majority of the initial contributions in TV-based models on manifolds were proposed in the early 2010s \cite{strekalovskiy2011total,cremers2013total,valkonen2013total}. These initial models had the major drawback that they were proposed for specific manifolds. In \cite{lellmann2013total}, the authors reformulated the variational problem on arbitrary Riemannian manifolds into a multi-label optimization problem.  In \cite{weinmann2014total}, the generalized ROF model was formulated in a fully intrinsic fashion. 

When entering the second half of the 2010s we see generalizations to several models emerging. A second-order model for cyclic data was proposed in \cite{bergmann2014second}, and soon, popular real-valued models were extended to the manifold case: a general second-order method \cite{bacak2016second}, infimal convolution models \cite{bergmann2017infimal, bergmann2018priors}, and TGV for manifold-valued imaging \cite{bredies2018total, bergmann2018priors}. At the same time specialized models were extended to applications such as inpainting \cite{bergmann2015inpainting}, segmentation \cite{weinmann2016mumford}, or manifold-valued inverse problems \cite{baust2016combined,storath2018wavelet}. Additionally, new problem settings arose with the emergence TV for manifold-valued data on graphs \cite{bergmann2018graph}.

As non-smooth models were extended to manifolds, so did the solvers. The proposed algorithms were typically proposed on \emph{Hadamard manifolds}: Riemannian manifolds that are complete, simply connected and have non-positive sectional curvature. On these manifolds, convexity can be generalized, which has been crucial in many of the convergence proofs. Proposed algorithms include the Proximal Point Algorithm (PPA) \cite{banert2014backward} and its extension, the Cyclic Proximal Point Algorithm (CPPA) \cite{bacak2014computing}; Iteratively Reweighted Least Squares (IRLS) \cite{grohs2014total,Bergmann2016} and an adaptation in \cite{grohs2016total}; the Parallel Douglas-Racheford Algorithm (PDRA) \cite{bergmann2016parallel}; and the recently derived exact and linearized Riemannian Chambolle Pock Algorithms (eRCPA/lRCPA) \cite{bergmann2019fenchel}. Despite the Hadamard constraint, it should be noted that these algorithms perform well on manifolds with positive curvature in experiments. Methods that do not rely on the Hadamard constraint include exact methods as in \cite{storath2016exact}, which uses exact solutions for $\ell^1$-TV with spherical data, and functional lifting \cite{vogt2019lifting}. These methods come with their own shortcomings: for the former this is limited applicability and for the latter this is high computational cost.


\paragraph{Non-smooth Optimization on Manifolds}
For non-smooth optimization on manifolds, the pioneering works include \cite{ferreira1998subgradient}, in which the subgradient method was extended to Riemannian manifolds and was shown to converge on Hadamard manifolds, and \cite{ferreira2002proximal}, whose authors extended the proximal map and proved convergence on Hadamard manifolds. Further development in non-smooth optimization was accelerated in the late 2000s by the introduction of the proximal subdifferential for manifolds \cite{azagra2005proximal},  the Clarke generalized subdifferential for manifolds \cite{hosseini2011generalized}, and a framework for duality on $CAT(0)$ metric spaces \cite{kakavandi2010duality}.


By the start of the 2010s, researchers started moving from Hadamard spaces to the non-Hadamard case. New developments include the introduction of a non-smooth version of the Kurdyka-Lojasiewicz (KL) inequality used to show the convergence of PPA on general Riemannian manifolds \cite{bento2011convergence}, convergence of subgradient descent for functions satisfying the KL inequality to a singular critical point \cite{hosseini2015convergence}, and a new approach to the convergence of PPA that extends previous results to a broader class of functions \cite{bento2016proximal}.

Around 2015, the first numerical implementations of these and new algorithms received attention: \cite{grohs2016nonsmooth} introduced a non-smooth trust region method for Riemannian manifolds and showed global convergence, \cite{grohs2016varepsilon} used an approximate subdifferential and proposed a descent method with global convergence, \cite{hosseini2017riemannian} proposed a gradient sampling algorithm and showed its global convergence, and \cite{hosseini2018line} proposed a line search algorithm and generalized the Wolfe conditions for Riemannian manifolds. Recently, a higher-order method was introduced: the Riemannian Semi-smooth Newton method (RSSN) \cite{de2018newton}, although to our knowledge no publicly-available implementation has existed so far. 
This method will form the basis of our approach.


\subsection{Contribution}
The contributions of this work are two-fold:

\paragraph{1. Primal-Dual Riemannian Semi-smooth Newton (PD-RSSN) for non-smooth optimization} While the RSSN method \cite{de2018newton} allows to find zeros of generic vector fields on manifolds in the same way that the Newton method allows to solve nonlinear systems of equations, the application to finding minimizers of non-smooth optimization problems is not straightforward, as the optimality conditions typically are in inclusion form. Using the recently proposed generalized Fenchel duality theory on manifolds \cite{bergmann2019fenchel}, we construct a primal-dual optimality system in the form of a vector field, which is then solved by RSSN. Overall, this provides a superlinearly convergent primal-dual scheme (PD-RSSN) for non-smooth Riemannian optimization.

\paragraph{2. Expanding the theoretical framework of Riemannian Semi-smooth Newton}
As in the classical Newton method, the proposed method requires to solve Newton-type system in each step. For larger-scale problems, solving for the Newton steps exactly is often not efficient. 
Therefore we propose an \emph{inexact} version of RSSN and provide a convergence proof (\cref{thm: inexact Newton}). We show that at least linear convergence can be expected in the inexact case. The theoretical results are validated by our numerical experiments \jla{on a dual (Huber-) regularized $\ell^2$-$TV$ denoising model on manifolds with positive and negative sectional curvature: The two-dimensional unit sphere $S^2$ and the manifold of symmetric positive semidefinite $3\times 3$ matrices $\mathcal{P}(3)$.}
The full source code for the reproducing the experiments in this work is available at
\vspace{.5em}\begin{center}
\href{https://github.com/wdiepeveen/Primal-Dual-RSSN}{\texttt{https://github.com/wdiepeveen/Primal-Dual-RSSN}}.
\end{center}\vspace{.5em}
The central concepts in this work, in particular the extension to the Riemannian setting, were developed by the first author in his thesis {\cite{diepeveen2020nonsmooth}}. Early considerations in the Euclidean case can be found in the thesis {\cite{rust2017}}.

%

\subsection{Outline}
In \cref{sec: prelim}, basic notation from differential geometry and Riemannian geometry is covered and the necessary definitions for manifold Fenchel duality theory are discussed along with the resulting non-smooth optimality systems. We also summarize the basic notions necessary for applying the Riemannian Semi-smooth Newton method. In \cref{sec: higher-order}, we discuss how RSSN can be used to solve the non-smooth optimality system from  \cref{sec: prelim}, which leads to the proposed PD-RSSN\ method. In \cref{sec: inexact}, we present the inexact Riemannian Semi-smooth Newton method and give a local convergence proof.  The application of our method to $\ell^2$-TV-like problems is discussed in \cref{sec: application}, and numerical results are shown in \cref{sec: numerics}. We conclude and summarize in \cref{sec:conclusions}.


\section{Preliminaries}
\label{sec: prelim}
In the first half of this section, the notation from differential and Riemannian geometry is summarized. In the second part, we briefly recapitulate manifold Fenchel duality theory and the Riemannian Semi-smooth Newton method.

%

\subsection{Notation}

For details regarding differential geometry and Riemannian geometry, we refer the reader to books such as \cite{lee2013smooth,sakai1996riemannian,carmo1992riemannian}. 

We write $\mathcal{M}$ and $\mathcal{N}$ for Riemannian manifolds. The tangent space at $p \in \mathcal{M}$ is denoted by $\mathcal{T}_p \mathcal{M}$ and for tangent vectors we write $X_p$ and $Y_p,$ or simply $X$ and $Y$. For the tangent bundle we have $\mathcal{T}\mathcal{M}:= \bigcup_{p\in \mathcal{M}} \mathcal{T}_p \mathcal{M}$. Similarly, we write $\mathcal{T}^*_p \mathcal{M}$ for \wdpa{the dual of the tangent space, or cotangent space}, $\xi_p$ and $\eta_p$ or simply $\xi$ and $\eta$ for covectors, and $\mathcal{T^*}\mathcal{M}:= \bigcup_{p\in \mathcal{M}} \mathcal{T}_p^* \mathcal{M}$ for the cotangent bundle. A cotangent vector $\xi \in \mathcal{T}_p^* \mathcal{M}$ acts on a tangent vector $X \in \mathcal{T}_p \mathcal{M}$ through the duality pairing $\langle\xi, X \rangle_p: = \xi(X) \in \mathbb{R}$. \wdpa{For a mapping $\Lambda: \mathcal{M} \rightarrow \mathcal{N}$, we write $D_p \Lambda: \mathcal{T}_p \mathcal{M} \rightarrow \mathcal{T}_{\Lambda(p)}\mathcal{N}$ for the differential of $\Lambda$ at $p\in \mathcal{M}$.}

We assume a Riemannian manifold is equipped with a metric. For some point $p \in \mathcal{M}$ the metric tensor
is denoted by $(\cdot, \cdot)_p: \mathcal{T}_p \mathcal{M} \times \mathcal{T}_p \mathcal{M} \rightarrow \mathbb{R}$ and the norm induced by the Riemannian metric is written as $\| \cdot \|_p$. The Riemannian distance between points $p,q\in \mathcal{M}$ is denoted by $d_{\mathcal{M}}(p,q)$. For the open metric ball of radius $r>0$ with center $p \in \mathcal{M}$ induced by this distance we write
\begin{equation}
        \mathcal{B}_{r}(p):=\left\{y \in \mathcal{M} \mid d_{\mathcal{M}}(p, q)<r\right\}.
\end{equation}
Through this metric, for $p \in \mathcal{M}$ we can define the \emph{musical isomorphisms} $\flat: \mathcal{T}_{p} \mathcal{M} \rightarrow \mathcal{T}_{p}^{*} \mathcal{M}$ as
\begin{equation}
        \langle X^{\flat}, Y\rangle_p=(X, Y)_{p} \text { for all } Y \in \mathcal{T}_{p} \mathcal{M}
\end{equation}
and its inverse $\sharp: \mathcal{T}_{p}^{*} \mathcal{M} \rightarrow \mathcal{T}_{p} \mathcal{M}$ as
\begin{equation}
        (\xi^{\sharp}, X)_{p}=\langle\xi, X\rangle_p \text { for all } X \in \mathcal{T}_{p} \mathcal{M}.
\end{equation}

The metric can also be used to construct a unique affine connection, the \emph{Levi-Civita connection} \jla{or \emph{covariant derivative,}} which is denoted by $\nabla_{(\cdot)}(\cdot)$. For geodesics $\gamma:[0,1] \rightarrow \mathcal{M,}$  we typically write  $\gamma_{p,q}(t)$ and $\gamma_{p,X}(t)$ to denote a minimizing geodesic connecting $p,q\in \mathcal{M}$ and a geodesic starting from $p\in\mathcal{M}$ with velocity $\dot{\gamma}_{p,X} (0) = X\in \mathcal{T}_p\mathcal{M}$. The latter defines the \emph{exponential map,} which is denoted by $\exp_p: \mathcal{G}_p  \rightarrow \mathcal{M}$ where $\mathcal{G}_p \subset \mathcal{T}_p\mathcal{M}$ is the set on which $\dot{\gamma}_{p,X}=:\exp_p(X)$ 
is defined. Furthermore, if $\mathcal{G}_p = \mathcal{T}_p\mathcal{M}$, the manifold is called \emph{complete.} For $\mathcal{G}'_p \subset \mathcal{G}_p$ a metric ball with radius $r_p$ on which $\exp_p$ is a diffeomorphism, the \emph{logarithmic map} is defined and we denote it by $\log_p: \exp(\mathcal{G}'_p ) \rightarrow \mathcal{G}'_p$. This radius $r_p$ is called the injectivity radius. Importantly, on simply connected, complete Riemannian manifolds with non-positive sectional curvature, the $\exp_p$ and $\log_p$ maps are globally defined. Such manifolds are called \emph{Hadamard manifolds.} Finally, we write \emph{parallel transport} of a vector $X \in \mathcal{T}_p\mathcal{M}$ from $p$ to $q \in \mathcal{M}$ as $\mathcal{P}_{q\leftarrow p}X$ and parallel transport of a covector $\xi \in \mathcal{T}_p^*\mathcal{M}$ as $\mathcal{P}_{q\leftarrow p}\xi$, where the latter is defined through the musical isomorphisms, as
\begin{equation}
        \mathcal{P}_{q\leftarrow p} \xi:=\left(\mathcal{P}_{q\leftarrow p} \xi^{\sharp}\right)^{\flat}.
\end{equation}


\subsection{Manifold Duality Theory}\label{sec: manifold duality}
In this section, generalizations of classical notions from non-smooth analysis to manifolds proposed in \cite{bergmann2019fenchel} are discussed. For the equivalent basic notions of convex analysis in vector spaces we refer the reader to \cite{rockafellar2009variational,clason2020introduction}. The section culminates in the recently derived primal-dual optimality conditions proposed by the authors of \cite{bergmann2019fenchel}. In the following, $\overline{\mathbb{R}}:=\mathbb{R} \cup\{\pm \infty\}$ denotes the extended real line.

\subsubsection{Non-smooth Analysis and a Fenchel Conjugation Scheme on Manifolds}
The notion of convexity can be defined on strongly convex subsets of Riemannian manifolds. 
\begin{definition}[strongly convex set, {\cite[Def. IV.5.1]{sakai1996riemannian}}]
         A subset $\mathcal{C} \subset \mathcal{M}$ of a Riemannian manifold $\mathcal{M}$ is said to be \emph{strongly convex} if, for all $p, q \in \mathcal{C}$, a minimal geodesic $\gamma_{p,q}$ between $p$ and $q$ exists, is unique, and lies completely in $\mathcal{C}$.
\end{definition}
The well-known notions of properness, convexity and lower semi-continuity can be generalized as follows.
\begin{definition}[proper, {\cite[Def. 2.11.i]{bergmann2019fenchel}}]
        A function $F : \mathcal{M}\rightarrow \overline{\mathbb{R}}$ is \emph{proper} if $\operatorname{dom} F := \{x \in \mathcal{M}| F(x) < \infty \} \neq \emptyset$ and $F(x) > -\infty$ holds for all $x \in \mathcal{M}$.
\end{definition}

\begin{definition}[convex, {\cite[Def. 2.11.ii]{bergmann2019fenchel}}]
        Suppose that $\mathcal{C} \subset \mathcal{M}$ is strongly convex. A proper function $F : \mathcal{M} \rightarrow \overline{\mathbb{R}}$ is called (geodesically) \emph{convex} on $\mathcal{C} \subset \mathcal{M}$ if, for all $p, q \in \mathcal{C,}$ the composition $F \circ \gamma_{p,q}(t)$ is a convex function on $[0, 1]$ in the classical sense.
\end{definition}

\begin{definition}[epigraph, {\cite[Def. 2.11.iii]{bergmann2019fenchel}}]
        Suppose that $\mathcal{A} \subset \mathcal{M}$. The \emph{epigraph} of a function $F : \mathcal{A} \rightarrow \overline{\mathbb{R}}$ is defined as
        \begin{equation}
        \operatorname{epi} F:=\{(x, \alpha) \in \mathcal{A} \times \mathbb{R} | F(x) \leq \alpha\}.
        \end{equation}
        
\end{definition}

\begin{definition}[lower semi-continuous, {\cite[Def. 2.11.iv]{bergmann2019fenchel}}]
        Suppose that $\mathcal{A} \subset \mathcal{M}$. A proper function $F : \mathcal{A} \rightarrow \overline{\mathbb{R}}$ is called \emph{lower semi-continuous} (lsc) if $\operatorname{epi} F$ is closed.
\end{definition}
The notions of subdifferentials and proximal mappings can also be extended to manifolds using the exponential map and the geodesic distance: 

\begin{definition}[subdifferential,{\cite[Def. 3.4.4]{udriste1994convex}}]
        Suppose that $\mathcal{C} \subset \mathcal{M}$ is strongly convex. The \emph{subdifferential} $\partial_{\mathcal{M}} F$ on $\mathcal{C}$ of a proper, convex function $F : \mathcal{C} \rightarrow \overline{\mathbb{R}}$  at a point $p\in  \mathcal{C}$ is defined as 
        \begin{equation}
        \partial_{\mathcal{M}} F(p):=\left\{\xi \in \mathcal{T}_{p}^{*} \mathcal{M} | F(q) \geq F(p)+\langle\xi, \log _{p} q\rangle_p \text { for all } q \in \mathcal{C}\right\}.
        \end{equation}
\end{definition}

\begin{definition}[proximal mapping, \cite{ferreira2002proximal}]
        Let $\mathcal{M}$ be a Riemannian manifold, $F : \mathcal{M} \rightarrow \overline{\mathbb{R}}$ be proper, and $\lambda > 0$. The \emph{proximal map} of $F$ is defined as 
        \begin{equation}
        \operatorname{prox}_{\lambda F}(p):=\arg \min_{q \in \mathcal{M}}\left\{\frac{1}{2\lambda } d_{\mathcal{M}}(p, q)^{2}+ F(q)\right\}.
        \end{equation}
\end{definition}
Tangent and cotangent spaces play an important role in the generalization of the Fenchel-dual functions. In this context, well-definedness of the exponential and the logarithmic map is ensured by restricting ourselves to the following subset, which is a localized variant of the pre-image of the exponential map:
\begin{definition}[{\cite[Def 2.8]{bergmann2019fenchel}}]
        Let $\mathcal{C} \subset \mathcal{M}$ and $p \in \mathcal{C}$. We define the tangent subset \linebreak $\mathcal{L}_{\mathcal{C}, p}\subset\mathcal{T}_{p}\mathcal{M}$ as
        \begin{equation}
                \mathcal{L}_{\mathcal{C}, p}:=\left\{X \in \mathcal{T}_{p} \mathcal{M} \mid \exp _{p} X \in \mathcal{C} \text { and }\|X\|_{p}=d_{\mathcal{M}}(\exp _{p} X, p)\right\}.
        \end{equation}
\end{definition}

With these basic notions, Fenchel duality theory can be generalized as in \cite{bergmann2019fenchel}. The Fenchel conjugate or Fenchel dual is defined by introducing a base point $m$ on the manifold.
\begin{definition}[$m$-Fenchel conjugate,{\cite[Def. 3.1]{bergmann2019fenchel}}]
        Suppose that $F : \mathcal{C} \rightarrow \overline{\mathbb{R}}$ and $m \in  \mathcal{C}$. The \emph{$m$-Fenchel conjugate} of $F$ is defined as the function $F^*_m
        : \mathcal{T}_p^*\mathcal{M} \rightarrow \overline{\mathbb{R}}$ such that
        \begin{equation}
        F_{m}^{*}\left(\xi_{m}\right):=\sup _{X \in \mathcal{L}_{C,m}}\left\{\left\langle\xi_{m}, X\right\rangle_m- F\left(\exp _{m} X\right)\right\}, \quad \xi_{m} \in \mathcal{T}_{m}^{*} \mathcal{M}.
        \end{equation}
\end{definition}
For the Fenchel biconjugate, we can then define the following.
\begin{definition}[$(mm')$-Fenchel biconjugate, {\cite[Def. 3.5]{bergmann2019fenchel}}]
                Suppose that $F : \mathcal{C} \rightarrow \overline{\mathbb{R}}$ and $m,m'\in  \mathcal{C}$. Then the \emph{$(mm')$-Fenchel biconjugate function} $F_{mm'}: \mathcal{C} \rightarrow \overline{\mathbb{R}}$ is defined as
                \begin{equation}
                F_{m m^{\prime}}^{* *}(p):=\sup _{\xi_{m^{\prime}} \in \mathcal{T}_{m^{\prime}}^{*} \mathcal{M}}\left\{\left\langle\xi_{m^{\prime}}, \log _{m^{\prime}} p\right\rangle_{m'}- F_{m}^{*}\left(\mathcal{P}_{m \leftarrow m^{\prime}} \xi_{m^{\prime}}\right)\right\}, \quad p \in \mathcal{C}.
                \end{equation}
\end{definition}
Basic properties and a more elaborate discussion of these generalized conjugate functions are covered in \cite{bergmann2019fenchel}. We will only focus on the relation between $F^{**}_{mm}$ and $F$. 

First, we note that geodesic convexity is often too strong a condition, so a weaker condition is used. Given a set $\mathcal{C} \subset \mathcal{M}$, $m \in \mathcal{C}$, and a function $F: \mathcal{C} \rightarrow \overline{\mathbb{R}},$ we define $f_{m}: \mathcal{T}_{m} \mathcal{M} \rightarrow \overline{\mathbb{R}}$ by
\begin{equation}
        f_{m}(X):=\left\{\begin{array}{cl}
                F\left(\exp _{m} X\right), & X \in \mathcal{L}_{\mathcal{C}, m}, \\
                +\infty, & X \notin \mathcal{L}_{\mathcal{C}, m}.
        \end{array}\right.
\end{equation}
It turns out that one can look at the convexity of the function $f_{m}: \mathcal{T}_{m} \mathcal{M} \rightarrow \overline{\mathbb{R}}$ in the usual  vector space sense on $\mathcal{T}_{m} \mathcal{M}$. For a more elaborate discussion of the discrepancy between the convexity of $f_m$ and (geodesic) convexity of $F$ we refer to \cite[Example 3.10]{bergmann2019fenchel}.

Now, the Fenchel-Moreau-Rockafellar theorem can also be extended to the manifold case:
\begin{theorem}[{\cite[Thm. 3.13]{bergmann2019fenchel}}]
        Suppose that $\mathcal{C} \subset \mathcal{M}$ is strongly convex and $m \in \mathcal{C}$. Let $F: \mathcal{C} \rightarrow \overline{\mathbb{R}}$ be proper. If $f_{m}$ is lsc and convex on $\mathcal{T}_{m} \mathcal{M},$ then $F=F_{m m}^{* *} .$ In this case $F_{m}^{*}$ is proper as well.
\end{theorem}

\subsubsection{The First-order Optimality Conditions}
\label{sec: optimality conditions}

In this section the first-order optimality conditions for a minimization problem of the form
\begin{equation}
        \inf_{p \in \mathcal{C}} \left\{ F(p)+G(\Lambda(p)) \right\}
        \label{eq: primal opt problem}
\end{equation}
are discussed. Here $\mathcal{C} \subset \mathcal{M}$ and $\mathcal{D} \subset \mathcal{N}$ are strongly convex sets, $F: \mathcal{C} \rightarrow \overline{\mathbb{R}}$ and $G: \mathcal{D} \rightarrow \overline{\mathbb{R}}$ are proper functions, and $\Lambda: \mathcal{M} \rightarrow \mathcal{N}$ is a differentiable map, possibly nonlinear, such that $\Lambda(\mathcal{C}) \subset \mathcal{D}$. Furthermore, we assume that $F: \mathcal{C} \rightarrow \overline{\mathbb{R}}$ is geodesically convex and that
\begin{equation}
        g_{n}(X):=\left\{\begin{array}{cl}
                G\left(\exp _{n} X\right), & X \in \mathcal{L}_{\mathcal{D}, n}, \\
                +\infty, & X \notin \mathcal{L}_{\mathcal{D}, n},
        \end{array}\right.
\end{equation}
is proper, convex and lsc on $\mathcal{T}_{n} \mathcal{N}$ for some $n \in \mathcal{D}$. 

Under these assumptions, \cref{eq: primal opt problem} can be rewritten into the following saddle-point formulation
\begin{equation}
        \inf_{p \in \mathcal{C}} \sup _{\xi_{n} \in \mathcal{T}_{n}^{*} \mathcal{N}} \left\{  \left\langle \xi_{n}, \log _{n} \Lambda(p)\right\rangle_n +F(p)-G_{n}^{*}\left(\xi_{n}\right) \right\}
\end{equation}

The authors of \cite{bergmann2019fenchel} propose two pairs of optimality conditions which we will also  refer to as the \emph{exact} and \emph{linearized} optimality conditions for the \emph{primal} and \emph{dual variables}. The proposed \emph{exact optimality conditions} are
\begin{align}
        \label{eq: exact primal optimality}
\mathcal{P}_{p\leftarrow m}\left(-(D_m \Lambda)^{*}\left[\mathcal{P}_{\Lambda(m) \leftarrow n} \xi_{n}\right]\right) & \in \partial_{\mathcal{M}} F(p), \\
\log _{n} \Lambda(p) & \in \partial G_{n}^{*}\left(\xi_{n}\right),
\end{align}
where $D_m \Lambda^{*}: \mathcal{T}_{\Lambda(m)}^{*} \mathcal{N} \rightarrow \mathcal{T}_{m}^{*} \mathcal{M}$ is the adjoint operator of $D_m \Lambda$. As shown in \cite{bergmann2019fenchel}, this system can be rewritten into~
\begin{align}
p&= \operatorname{prox}_{\sigma F} \left( \exp _{p}\left(\mathcal{P}_{p\leftarrow m}\left(-\sigma(D_m \Lambda)^{*}\left[\mathcal{P}_{ \Lambda(m) \leftarrow n} \xi_{n}\right]\right)^{\sharp}\right)\right),\\
\xi_{n}&= \operatorname{prox}_{\tau G^*_n} \left(\xi_{n}+\tau\left(\log _{n} \Lambda\left(p\right)\right)^{\flat} \right).\label{eq: exact update xin}
\end{align}
As coined in \cite{valkonen2014primal,bergmann2019fenchel}, \jla{we use the term} \enquote{exact} \jla{to refer to the fact that the operator $\Lambda$ is used in its original form in without linearization in} \cref{eq: exact update xin}. \jla{However, note that the linearization}
\begin{equation}
        \Lambda(p) \approx \exp _{\Lambda(m)} D \Lambda(m)\left[\log _{m} p\right]
        \label{eq: linearization}
\end{equation}
\jla{was still needed to obtain {\cref{eq: exact primal optimality}}, in particular so that the adjoint operator $D_m \Lambda$ can be constructed.}

The \emph{linearized optimality conditions} can be obtained by linearizing both the primal and the dual optimality condition. That is, for $n:= \Lambda(m),$ we want to solve
\begin{equation}
        \inf_{p \in \mathcal{M}} \inf_{\xi_{n} \in \mathcal{T}_{n}^* \mathcal{N}} \{F(p) + \left\langle \xi_{n}, D_m \Lambda\left[\log _{m} p\right]\right\rangle_n - G_{n}^{*}\left(\xi_{n}\right)\},
\end{equation}
and obtain the optimality system for the general $n\in\mathcal{N}$ case \cite{bergmann2019fenchel}
\begin{align}
\mathcal{P}_{p \leftarrow m}\left(-(D_m \Lambda)^{*}\left[\mathcal{P}_{\Lambda(m) \leftarrow n} \xi_{n}\right]\right) & \in \partial_{\mathcal{M}} F(p), \\
D_m \Lambda\left[\log _{m} p\right] & \in \partial G_{n}^{*}\left(\xi_{n}\right),
\end{align}
which can be rewritten into
\begin{align}
\label{eq: lin primal prox}
p &= \operatorname{prox}_{\sigma F}\left(\exp _{p}\left(\mathcal{P}_{ p \leftarrow m} \left(-\sigma \left(D_m \Lambda\right)^{*}\left[\mathcal{P}_{\Lambda(m) \leftarrow n} \xi_{n}\right]\right)^{\sharp}\right)\right), \\
\xi_{n} &= \operatorname{prox}_{\tau G_{n}^{*}}\left(\xi_{n}+\tau\left(\mathcal{P}_{n \leftarrow \Lambda(m)} D_m \Lambda\left[\log _{m} p\right]\right)^{\flat}\right).
\label{eq: lin dual prox}
\end{align}
\jla{We note that  the connection between the linearized optimality conditions and a solution of the original problem has only been conclusively established in the Euclidean case so far {\cite{valkonen2014primal}.} In the Riemannian case, the following weak duality result holds:} 
\begin{theorem}[{\cite[Thm. 4.2]{bergmann2019fenchel}}]
        Let $n:= \Lambda(m)$. The dual problem of 
                \begin{equation}
                \inf_{p \in \mathcal{M}} \left\{F(p)+G(\exp _{\Lambda(m)} D_m \Lambda\left[\log _{m} p\right]) \right\}
        \end{equation}
        is given by
        \begin{equation}
                \sup _{\xi_{n} \in \mathcal{T}_{n}^* \mathcal{N}} \left\{ F_{m}^{*}\left(-(D_m \Lambda)^{*}\left[\xi_{n}\right]\right)-G_{n}^{*}\left(\xi_{n}\right) \right\}
        \end{equation}
        and weak duality holds, i.e.,
        \begin{equation}
        \inf_{p \in \mathcal{M}} \left\{ F(p)+G(\exp _{\Lambda(m)} D_m \Lambda\left[\log _{m} p\right]) \right\} \geq \sup_{\xi_{n} \in \mathcal{T}^*_{n} \mathcal{N}} \left\{-F_{m}^{*}\left(-(D_m \Lambda)^{*}\left[\xi_{n}\right]\right)-G_{n}^{*}\left(\xi_{n}\right) \right\}.
        \end{equation}
\end{theorem}
If one performs a suitable fixed-point iteration on \cref{eq: lin primal prox,eq: lin dual prox}, one obtains the \emph{linearized Riemannian Chambolle Pock algorithm,} which has been shown to converge on Hadamard manifolds \cite[Thm. 4.3]{bergmann2019fenchel}. In this work, we aim to construct a higher-order method instead.

\subsection{The Riemannian Semi-smooth Newton Method}
\label{sec: RSSN method}
This section is devoted to the \emph{Riemannian Semi-smooth Newton} (RSSN) method, which we propose to use for solving \cref{eq: primal opt problem}. 

\begin{remark}
        \label{rem: M rewrite}
        Throughout this section and \cref{sec: inexact}, we will develop RSSN on some manifold $\mathcal{M}'$. This is not the original manifold $\mathcal{M}$ defining the feasible set;\  we will later apply the RSSN method to $\mathcal{M}' = \mathcal{M} \times \mathcal{T}^*_{n} \mathcal{N}$ (and $p' = (p,\xi_n)$). In order not to clutter the notation, we still write $\mathcal{M}$ instead of $\mathcal{M}'$ (and $p$ instead of $p'$) throughout this section and  \cref{sec: inexact}.
\end{remark}

In the real-valued case, the Newton method is formulated for solving a smooth system of equations $X(p) = 0$ for some non-linear map $X:\mathbb{R}^d \rightarrow \mathbb{R}^d$, where -- in the optimization context -- this map implements optimality conditions for the minimization problem. \wdpa{In the real-valued case, the classical Newton iteration is given by}
\begin{equation}
        p^{k+1} = p^{k} - \nabla X(p^k)^{-1} X(p^k).
\end{equation}

In the manifold case, there is no unique generalization. One possibility is to consider a \emph{vector field} $X: \mathcal{M} \rightarrow \mathcal{T} \mathcal{M}$ and define the Newton iteration through the covariant derivative:
\begin{equation}
p^{k+1} = \exp_{p^k}(-[\nabla_{(\cdot)} X]_{p^k}^{-1}X(p^k)) 
\end{equation}
in order to find a zero $0 = X(p)$ of the vector field. In this section, we will focus on a generalized covariant derivative approach for finding zeros of \emph{semi-smooth} vector fields based on \cref{eq: lin primal prox}--\cref{eq: lin dual prox}. Throughout this section we will use the notions and results developed in~\cite{de2018newton}.

\begin{definition}[(locally) Lipschitz, {\cite[Def. 6]{de2018newton}}]
        A vector field $X$ on $\mathcal{M}$ is said to be \emph{Lipschitz continuous} on $\Omega \subset \mathcal{M}$ if there exists a constant $L>0$ such that, for all $p, q \in \Omega$, there holds
        \begin{equation}
                \left\|\mathcal{P}_{q \leftarrow p} X(p)-X(q)\right\|_q \leq L d_{\mathcal{M}}(p,q), \quad \forall p, q \in \Omega,
        \end{equation}
        Moreover, for a point $p \in \mathcal{M}$, if there exists $\delta>0$ such that $X$ is Lipschitz continuous on the open ball $B_{\delta}(p),$ then $X$ is said to be Lipschitz continuous at p. Moreover, if for all $p \in \mathcal{M}$, $ X$ is Lipschitz continuous at $p$, then $X$ is said to be \emph{locally Lipschitz continuous} on $\mathcal{M}$.
\end{definition}
Subsequently, we can generalize Rademacher's theorem to Lipschitz vector fields.
\begin{theorem}[{\cite[Thm. 10]{de2018newton}}]
        If $X$ is a locally Lipschitz continuous vector field on $\mathcal{M}$, then $X$ is almost everywhere differentiable on $\mathcal{M}$.
\end{theorem}
Hence, it makes sense to define the generalized covariant derivative. 
\begin{definition}[Clarke generalized covariant derivative, {\cite[Def. 11]{de2018newton}}]
         The \emph{Clarke generalized covariant derivative} $\partial_{C} X$ of a locally Lipschitz continuous vector field $X$ at a point $p \in \mathcal{M}$ is defined as the set-valued mapping $\partial_{ C} X\left(p\right):\mathcal{T}_{p} \mathcal{M} \rightrightarrows \mathcal{T}_{p} \mathcal{M}$,
         \begin{equation}
         \partial_{C} X\left(p\right):=\operatorname{co}\left\{V \in \mathcal{L}\left(\mathcal{T}_{p} \mathcal{M}\right): \exists (p^{k})_{k\geq 0} \subset \mathcal{D}_{X}, \lim _{k \rightarrow\infty} p^{k}=p, \; V=\lim _{k \rightarrow\infty} \mathcal{P}_{p \leftarrow p^{k} } \nabla X(p^{k})\right\},
         \end{equation}
        where $\mathcal{D}_{X} \subset \mathcal{M}$ is the set on which $X$ is differentiable, \enquote{co} represents the convex hull and $\mathcal{L}\left(\mathcal{T}_{p} \mathcal{M}\right)$ is the vector space consisting of all bounded linear operators from $\mathcal{T}_{p}\mathcal{M}$ to $\mathcal{T}_{p} \mathcal{M}$.
\end{definition}
In addition to local Lipschitzness, we will also need the \emph{directional derivative} for the notion of semi-smoothness. We follow the definition in \cite{de2018newton}.
\begin{definition}[directional derivative]
        The \emph{directional derivative} of a vector field $X$ on $\mathcal{M}$ at $p \in \mathcal{M}$ in the direction $v \in \mathcal{T}_{p} \mathcal{M}$ is defined by
        \begin{equation}
                X'(p, v):=\lim _{t \searrow 0} \frac{1}{t}\left[\mathcal{P}_{p \leftarrow \exp _{p}(t v)} X\left(\exp _{p}(t v)\right)-X(p)\right] \in \mathcal{T}_{p} \mathcal{M},
        \end{equation}
        whenever the limit exists. If this directional derivative exists for every $v$, then $X$ is said to be \emph{directionally differentiable} at $p$.
\end{definition}

Finally, we are able to generalize the notion of semi-smoothness to vector fields. 
\begin{definition}[semi-smooth vector field, {\cite[Def. 18]{de2018newton}}]
        \label{def: semismooth manifold}
        A vector field $X$ on $\mathcal{M}$ that is Lipschitz continuous at $p \in \mathcal{M}$ and directionally differentiable at $q \in B_{\delta}\left(p\right)$ for all directions in $\mathcal{T}_{p} \mathcal{M}$, is said to be \emph{semi-smooth} at $p$ iff for every $\epsilon>0$ there exists $0<\delta<r_{p}$, where $r_{p}$ is the injectivity radius, such that
        \begin{equation}
                \|X\left(p\right)-\mathcal{P}_{ p \leftarrow q}\left[X(q)+V_{q} \log_{q} p\right]\|_p \leq \epsilon d_{\mathcal{M}}\left(p, q\right), \quad \forall q \in B_{\delta}\left(p\right), \quad \forall V_{q} \in \partial_{C} X(q).
                \label{eq: semismoothness}
        \end{equation}
        The vector field $X$ is said to be \emph{$\mu$-order semi-smooth} at $p$ for $0<\mu \leq 1$ iff there exist $\epsilon>0$ and $0<\delta<r_{p}$ such that
        \begin{equation}
                \|X\left(p\right)-\mathcal{P}_{p \leftarrow q }\left[X(q)+V_{q} \log_{q} p\right]\|_p \leq \epsilon d_{\mathcal{M}}\left(p, q\right)^{1+\mu}, \quad \forall q \in B_{\delta}\left(p\right), \quad \forall V_{q} \in \partial_{C}X(q).
        \end{equation}
\end{definition}
        

The Riemannian Semi-smooth Newton method for finding a zero of a vector field, i.e., $X(p) = 0$, is shown in \cref{alg: RSSN}. It extends on the classical Newton method by replacing the  classical Jacobian $\nabla X(p^k)$ by an element from the Clarke generalized covariant derivative $\partial_C X(p^k)$, and performing the update step using the exponential map.

\begin{algorithm}[h!]
        \caption{Riemannian Semi-smooth Newton}
        \label{alg: RSSN}
        \begin{algorithmic}
                \STATE{Initialization: $p^0 \in \mathcal{M}, k := 0$}
                \WHILE{not converged}
                \STATE{Choose any $V(p^k) \in \partial_{C} X(p^k)$}
                \STATE{Solve $V(p^k) d^k= -X(p^k)$ in the vector space $\mathcal{T}_{p^k}\mathcal{M}$}
                \STATE{$p^{k+1} := \exp_{p^k}(d^k)$}
                \STATE{$k := k+1$}
                \ENDWHILE
        \end{algorithmic}
\end{algorithm}

We have the following local convergence result.
\begin{theorem}[{\cite[Thm. 19]{de2018newton}}]
        \label{thm: RSSN conv}
        Let $X$ be a locally Lipschitz-continuous vector field on $\mathcal{M}$ and $p^{*} \in \mathcal{M}$ be a solution of the problem $X(p)=0$. Assume that $X$ is semi-smooth at $p^*$ and that all $V_{p^{*}} \in \partial_{C} X(p^{*})$ are invertible. Then there exists a $\delta>0$ such that, for every starting point $p^{0} \in B_{\delta}\left(p^{*}\right) \backslash\left\{p^{*}\right\}$, the sequence  $(p^{k})_{k\geq 0}$ generated by \cref{alg: RSSN} is well-defined, contained in $B_{\delta}\left(p^{*}\right)$ and converges superlinearly to $p^{*}$. If additionally $X$ is $\mu$-order semi-smooth at $p^{*},$ then the convergence of $(p^{k})_{k\geq 0}$ to $p^{*}$ is of order $1+\mu$.
\end{theorem}

\section{A Higher-order Primal-dual Method for Manifolds (PD-RSSN)}
\label{sec: higher-order}

%
%
We will now merge the ideas of Riemannian Semi-smooth Newton (RSSN) and Fenchel duality theory on manifolds in order to solve the original problem 
\begin{equation}
\inf_{p \in \mathcal{M}} \{F(p) +  G(\Lambda(p))\}.
\label{eq: exact problem}
\end{equation}
Using linearization, we saw in \cref{sec: manifold duality} that solving \cref{eq: exact problem} can be approximated by solving 
\begin{equation}
        \inf_{p \in \mathcal{M}} \sup_{\xi_{n} \in \mathcal{T}_{n}^* \mathcal{N}} \{F(p) + \left\langle \xi_{n}, D_m \Lambda\left[\log _{m} p\right] \right\rangle_n - G_{n}^{*}\left(\xi_{n}\right)\},
\end{equation}
which is characterized by the optimality conditions
\begin{align}
        \label{eq: primal eq}
        p &= \operatorname{prox}_{\sigma F}\left(\exp _{p}\left(\mathcal{P}_{p \leftarrow m} \left(-\sigma \left(D_m \Lambda\right)^{*}\left[\mathcal{P}_{\Lambda(m) \leftarrow n} \xi_{n}\right]\right)^{\sharp}\right)\right), \\
        \label{eq: dual eq}
        \xi_{n} &= \operatorname{prox}_{\tau G_{n}^{*}}\left(\xi_{n}+\tau\left(\mathcal{P}_{n \leftarrow \Lambda(m) } D_m \Lambda\left[\log _{m} p\right]\right)^{\flat}\right),
\end{align}
where $\sigma,\tau > 0$. In this section we focus on rewriting this system of non-linear equations into a form that is amenable to the Riemannian Semi-smooth Newton method, i.e., into the problem of finding a zero of a vector field. 

\label{sec: choose Newton}
While the dual variable $\xi_n$ lives in a vector space and \eqref{eq: dual eq} immediately translates into 
\begin{equation}
\xi_{n} - \operatorname{prox}_{\tau G_{n}^{*}}\left(\xi_{n}+\tau\left(\mathcal{P}_{n \leftarrow \Lambda(m) } D_m \Lambda\left[\log _{m} p\right]\right)^{\flat}\right) = 0,
\end{equation}
the primal optimality condition \eqref{eq: primal eq} is an equation in $\mathcal{M}$. In order to obtain a vector field form, we apply the logarithm, which yields elements from the tangent bundle\ $\mathcal{T}\mathcal{M}$: 
\begin{equation}
\label{eq: VF approach}
        -\log_p \operatorname{prox}_{\sigma F}\left(\exp _{p}\left(\mathcal{P}_{p \leftarrow m} \left(-\sigma \left(D_m \Lambda\right)^{*}\left[\mathcal{P}_{\Lambda\left(m\right) \leftarrow n} \xi_{n}\right]\right)^{\sharp}\right)\right) = 0.
\end{equation}
Such a zero implies \cref{eq: primal eq}. Then, we define the vector field $X: \mathcal{M} \times  \mathcal{T}_{n}^* \mathcal{N} \rightarrow \mathcal{T}\mathcal{M} \times \mathcal{T}_{n}^* \mathcal{N}$  as
\begin{equation}
X(p,\xi_n) := \begin{pmatrix}
-\log_p \operatorname{prox}_{\sigma F}\left(\exp _{p}\left(\mathcal{P}_{p \leftarrow m} \left(-\sigma \left(D_m \Lambda\right)^{*}\left[\mathcal{P}_{\Lambda\left(m\right) \leftarrow n} \xi_{n}\right]\right)^{\sharp}\right)\right)\\ 
\xi_{n} - \operatorname{prox}_{\tau G_{n}^{*}}\left(\xi_{n}+\tau\left(\mathcal{P}_{n \leftarrow \Lambda\left(m\right)} D_m \Lambda\left[\log _{m} p\right]\right)^{\flat}\right)
\end{pmatrix}.
\label{eq: X vf}
\end{equation}
which allows to apply RSSN. We refer to this approach as \emph{Primal-Dual Riemannian Semi-smooth Newton} (PD-RSSN).  


\begin{remark}
        Here it becomes obvious why first transforming the set-valued equivalent optimality conditions in \cref{eq: lin primal prox,eq: lin dual prox} into the prox-based equality form is crucial: we can now \emph{differentiate} $X$ in order to apply the Newton-based RSSN\ method.
\end{remark}

If the vector field $X$ in \cref{eq: X vf} is smooth, we obtain the covariant derivative \cite[Sec. 6.4.2]{diepeveen2020nonsmooth} \wdpa{at $(Y_p, \eta_{\xi_n}) \in \mathcal{T}_p\mathcal{M}\times \mathcal{T}_n^*\mathcal{N}$}:
\begin{align}
        \nabla_{(Y_p, \eta_{\xi_n})} X    = \begin{pmatrix}
                -\nabla_{Y_{p}}\log_{(\cdot)} f_1(p,\xi_n)  - D_p \log_{p} f_1(\cdot,\xi_n)[Y_p] - D_{\xi_n} \log_{p} f_1(p,\cdot) [\eta_{\xi_n}] \nonumber \\
         - D_p f_2(\cdot,\xi_n)[Y_p]+ \eta_{\xi_n} -  \nabla_{\eta_{\xi_n}} f_2(p,\cdot)
    \end{pmatrix},
\end{align}
where
\begin{align}
        f_1(p,\xi_n) & := \operatorname{prox}_{\sigma F}\left(\exp _{p}\left(\mathcal{P}_{p \leftarrow m} \left(-\sigma \left(D_m \Lambda\right)^{*}\left[\mathcal{P}_{ \Lambda\left(m\right) \leftarrow n} \xi_{n}\right]\right)^{\sharp}\right)\right),\\
        f_2(p,\xi_n) & := \operatorname{prox}_{\tau G_{n}^{*}}\left(\xi_{n}+\tau\left(\mathcal{P}_{n \leftarrow      \Lambda\left(m\right) } D_m \Lambda\left[\log _{m} p\right]\right)^{\flat}\right).
\end{align}
However, we generally cannot assume smoothness of $X$ and therefore require  a generalization of the differential. First, consider the following generalization of Rademacher's theorem.

\begin{theorem}
        Let $\mathcal{M}$ and $\mathcal{N}$ be smooth manifolds. If $F:\mathcal{M}\rightarrow \mathcal{N}$ is a locally Lipschitz continuous function, then $F$ is almost everywhere differentiable on $\mathcal{M}$.
\end{theorem}
\begin{proof}
        The proof is the same as  \cite[Thm. 10]{de2018newton} with a general manifold $\mathcal{N}$ instead of $\mathcal{T}\mathcal{M}$.
\end{proof}
The generalized differential can now be defined as follows; see \cite{hosseini2011generalized} for a related definition in the more restricted case of functions mapping from a manifold into the real values:
\begin{definition}[Clarke generalized differential]
        The \emph{Clarke generalized differential} $D_C F$ of a locally Lipschitz continuous function $F:\mathcal{M}\rightarrow \mathcal{N}$ at $p \in \mathcal{M}$ is defined as the set-valued mapping $D_{C} F(p):\mathcal{T}_{p} \mathcal{M}\rightrightarrows \mathcal{T}_{F(p)} \mathcal{N}$,
        \begin{equation}
        D_{C} F(p):=\operatorname{co}\left\{V \in \mathcal{L}\left(\mathcal{T}_{p} \mathcal{M},\mathcal{T}_{F(p)} \mathcal{N}\right): \exists (p^{k})_{k\geq 0} \subset \mathcal{D}_{F}, \lim _{k \rightarrow\infty} p^{k}=p, \; V=\lim _{k \rightarrow\infty} D_{p^k} F \right\},
        \end{equation}
        where $\mathcal{D}_{F}\subset \mathcal{M}$ is the set on which $F$ is differentiable, \enquote{co} represents the convex hull and $\mathcal{L}\left(\mathcal{T}_{p} \mathcal{M},\mathcal{T}_{F(p)} \mathcal{N}\right)$ denotes the vector space consisting of all bounded linear operator from $\mathcal{T}_{p}\mathcal{M}$ to $\mathcal{T}_{F(p)} \mathcal{N}$.
\end{definition}


In block notation, the generalized covariant derivative of $X$ is of the form

\begin{equation}
        \partial_{C} X(p,\xi_n) = \begin{bmatrix}
                -\partial_{C} (\log_{(\cdot)} f_1(p,\xi_n) ) (p)- D_C (\log_{p} f_1(\cdot,\xi_n)) (p)& -D_{C} (\log_{p} f_1(p,\cdot)) (\xi_n)\\ 
                -D_{C} (f_2(\cdot,\xi_n) )(p) & I -\partial_{C} (f_2(p,\cdot)) (\xi_n)
        \end{bmatrix}.
        \label{eq: Newton}
\end{equation}
Under the special assumption that the manifold of interest is \emph{symmetric}, the four components in \cref{eq: Newton} can be computed explicitly up to the differentials of the proximal maps using Jacobi fields \cite[Lemma 2.3]{persch2018optimization}. For the technical details we refer to the first author's thesis \cite[Sec. 6.4.3 and 6.4.4]{diepeveen2020nonsmooth}. For typical imaging applications this is convenient, as many typical manifolds of interest, in particular $S^n$ and $\mathcal{P}(n)$, are symmetric. 

\begin{remark}
        Semi-smoothness of $X$ does not follow directly from the proposed construction. In general it is problem-specific but relatively straightforward to prove, as it can be shown that piecewise smooth functions are semi-smooth. 
\end{remark}


\section{Inexact Riemannian Semi-smooth Newton}
\label{sec: inexact}
Before moving to applications, we will consider a generalization of RSSN, the \emph{Inexact Riemannian Semi-smooth Newton} (IRSSN) method in \cref{alg: inexact RSSN}. In contrast to the exact version, it only requires to solve the linear system up to a residual term $r^k$. The reader should bear in mind that \cref{rem: M rewrite} still holds for the remainder of this section. The main motivation for this method is that solving the Newton system with high precision can be very expensive. Solving the system inexactly, for example using an iterative method, can potentially ameliorate this problem. 

\begin{algorithm}[H]
        \caption{Inexact Riemannian Semi-smooth Newton}
        \label{alg: inexact RSSN}
        \begin{algorithmic}
                \STATE{Initialization: $p^0 \in \mathcal{M},a^{0} \geq 0, k := 0$}
                \WHILE{not converged}
                \STATE{Choose $V_k(p^k) \in \partial_{C} X(p^k)$}
                \STATE{Solve $V_k(p^k) d^k= -X(p^k) + r^k$ in $\mathcal{T}_{p^k}\mathcal{M}$ where $\|r^k\|_{(p^k)} \leq a^{k} \|X(p^{k})\|_{(p^k)}$}
                \STATE{$p^{k+1} := \exp_{p^k} (d^k)$}
                \STATE{Choose $a^{k+1} \geq 0$}
                \STATE{$k := k+1$}
                \ENDWHILE
        \end{algorithmic}
\end{algorithm}

In this section, we focus on proving \cref{thm: inexact Newton}: a local convergence result for \cref{alg: inexact RSSN} on Riemannian manifolds. The proof is based on the ideas of the Inexact Semi-smooth Newton methods in $\mathbb{R}^n$ as discussed in \cite{martinez1995inexact,facchinei1996inexact}. The technical details are inspired by the convergence proof for Riemannian Semi-smooth Newton \cite{de2018newton}.

\subsection{Towards a Convergence Proof for Inexact Riemannian Semi-smooth Newton}

Starting with the technicalities, we first need to account for curvature. In particular, we need to account for how geodesics spread, which we can formalize in a single quantity \cite[Def. 2]{de2018newton}:
\begin{definition}[{\cite[Def. 2]{de2018newton}}]
        \label{def: Kp}
        Let $p \in \mathcal{M}$ and $r_{p}$ be the radius of injectivity of $\mathcal{M}$ at $p .$ Define the quantity
        \begin{equation}
                K_{p}:=\sup \left\{\frac{d_{\mathcal{M}}\left(\exp _{q} u, \exp _{q} v\right)}{\|u-v\|_q}: q \in B_{r_{p}}(p), u, v \in \mathcal{T}_{q} \mathcal{M}, u \neq v,\|v\|_q \leq r_{p},\|u-v\|_q \leq r_{p}\right\}.
        \end{equation}
\end{definition}
The following remark from \cite[Remark 3]{de2018newton} provides some intuition:
\begin{remark}
        \label{rem: Kp}
This number $K_{p}$ measures how fast the geodesics spread apart in $\mathcal{M}$. In particular, when $u=0 \in \mathcal{T}_q \mathcal{M}$ or more generally when $u$ and $v$ are on the same line through $0$, $d_{\mathcal{M}}\left(\exp _{q} u, \exp _{q} v\right)=$ $\|u-v\|_q$. Hence, $K_{p} \geq 1,$ for all $p \in \mathcal{M}$. When $\mathcal{M}$ has non-negative sectional curvature, the geodesics spread apart less than the rays, i.e., $d_{\mathcal{M}}\left(\exp _{p} u, \exp _{p} v\right) \leq\|u-v\|_q$ and, in this case, $K_{p}=1$ for all $p \in \mathcal{M}$.
\end{remark}

Next, remember the definition of an operator norm.
\begin{definition}[{\cite[Def. 4]{de2018newton}}]
        Let $p \in \mathcal{M}$. The \emph{norm of a linear map} $A: \mathcal{T}_{p} \mathcal{M} \rightarrow \mathcal{T}_{p} \mathcal{M}$ is defined by
        \begin{equation}
                \|A\|_p:=\sup \left\{\|A v\|_p: v \in  \mathcal{T}_{p} \mathcal{M},\|v\|_p\leq 1\right\}.
        \end{equation}
\end{definition}

We have the following result. 
\begin{lemma}[{\cite[Lemma 17]{de2018newton}}]
        \label{lemma: V bound}
        Let $X$ be a locally Lipschitz continuous vector field on~$\mathcal{M}$. Assume that all elements $V_{p} \in \partial_{C} X(p)$ are invertible at base point $p \in \mathcal{M}$ and let \hfill \linebreak 
        $\lambda_{p} \geq \max \left\{\|V_{p}^{-1}\|_p: V_{p} \in \partial_{C} X(p)\right\}$. Then, for every $\epsilon>0$ satisfying $\epsilon \lambda_{p}<1$, there exists $0<\delta<r_{p}$ such that all $V_{q} \in \partial_{C} X(q)$ are invertible on $B_{\delta}\left(p\right)$ and
        
        \begin{equation}
        \label{eq: inv bdd}
        \|V_{q}^{-1}\|_q \leq \frac{\lambda_{p}}{1-\epsilon \lambda_{p}}, \quad \forall q \in B_{\delta}\left(p\right), \quad \forall V_{q} \in \partial_{C} X(q).
        \end{equation}
        
\end{lemma}

\subsection{Fast Local Convergence for Semi-smooth Vector Fields}
With these tools we can move on to the main result of this section.
\begin{theorem}
        \label{thm: inexact Newton}
        Let $X$ be locally Lipschitz continuous vector field on $\mathcal{M}$ and $p^{*} \in \mathcal{M}$ be a solution of problem the $X(p)=0$. Assume that $X$ is semi-smooth at $p^{*}$ and that all $V_{p^{*}} \in \partial_{C} X(p^{*})$ are invertible. Then the following statements hold:
        \begin{enumerate}[label=(\roman*)]
                \item There exist $a >0$ and $\delta>0$ such that, for every $p^{0} \in B_{\delta}\left(p^{*}\right)$ and $a^k \leq a$, the sequence $(p^{k})_{k\geq 0}$ generated by \cref{alg: inexact RSSN} is well-defined, is contained in $B_{\delta}\left(p^{*}\right)$ and converges Q-linearly to the solution $p^{*}$.
                \item If the sequence $(p^{k})_{k\geq 0}$ generated by \cref{alg: inexact RSSN} converges to the solution $p^{*}$ and further $\|r^{k}\|_{(p^k)}\in o\left(\|X(p^{k})\|_{(p^k)}\right)$, then the rate of convergence is $Q$-superlinear.
                \item If the sequence $(p^{k})_{k\geq 0}$ generated by \cref{alg: inexact RSSN} converges to the solution $p^{*}$, $X$ is $\mu$-order semi-smooth at $p^{*}$, and $\|r^{k}\|_{(p^k)}\in O\left(\|X(p^{k})\|_{(p^k)}^{1+\mu}\right)$, then the rate of convergence is of $Q$-order $1+\mu$.
        \end{enumerate}
\end{theorem}

\begin{proof}
        (i) 
        Let $K_{p^*}$ be as defined in \cref{def: Kp} and let $r_{p^*}$ be the injectivity radius. Since $X$ is locally Lipschitz, there exist constants $\hat{\delta} > 0$ and $L$ such that, for all $p \in B_{\hat{\delta}}(p^*),$
        \begin{equation}
        \label{eq: thm - lipschitz}
        \|X(p)\|_p = \|\mathcal{P}_{p \leftarrow p^*}X(p^*) - X(p)\|_{p} \leq L d_{\mathcal{M}}(p,p^*).
        \end{equation}
        The equality holds since $X(p^*) = 0$ and parallel transport is linear.
        
        Now, since all $V_{p^{*}} \in \partial_{C} X(p^{*})$ are invertible at $p^* \in \mathcal{M}$ by assumption, we can take $\lambda_{p^*} \geq \max \{ \| V_{p^*}^{-1}\|_{p^*}: V_{p^*} \in \partial_{C} X(p^*) \}$. Furthermore, take $a < \frac{1}{\lambda_{p^{*}} L K_{p^{*}} }$, choose $a^k\leq a $ $\forall k \in \mathbb{N}$ and $\epsilon$ satisfying $\epsilon \lambda_{p^{*}}\left(1+K_{p^{*}}\right)<1 - a \lambda_{p^{*}} L K_{p^{*}}$. As $\epsilon \lambda_{p^{*}} <1$, by \cref{lemma: V bound} we can find a $0 < \delta < \min\{\hat{\delta},r_{p^*}\}$ such that, for all $p \in B_{\delta}(p^*)$ and $V_{p} \in \partial_{C} X(p),$ 
        \begin{align}
        \label{eq: V bound}
        \|V_{p}^{-1}\|_p  \leq \frac{\lambda_{p^*}}{1-\epsilon \lambda_{p^*}}.
        \end{align}
        From the semi-smoothness of $X$,
        \begin{align}
        \label{eq: semismooth}
        \|X(p^*)-\mathcal{P}_{p^* \leftarrow p}\left[X(p)+V_{p} \log_{p} p^*\right]\|_{(p^*)} \leq \epsilon d_{\mathcal{M}}(p, p^*)
        \end{align}
        holds according to \cref{eq: semismoothness}.
        
        We now show that for our chosen $\delta$ the Newton iteration is well-defined. Let $k \in \mathbb{N}$ and assume that $p^k \in B_{\delta}(p^*)$. Let $d^k$ be such that
        \begin{equation}
        \|V_{p^k} d^k + X(p^k)\|_{(p^k)} \leq a^k\|X(p^k)\|_{(p^k)}.
        \label{eq: residual ineq}
        \end{equation}
        Then
        \begin{align}
        \| \log_{p^k} p^* - d^k \|_{(p^k)} & = \| \log_{p^k} p^* + V_{p^k}^{-1} X(p^k) - V_{p^k}^{-1}(V_{p^k} d^k + X(p^k) ) \|_{(p^k)} \\
        & \leq \|\log_{p^k} p^* + V_{p^k}^{-1} X(p^k)\|_{(p^k)} + \|V_{p^k}^{-1}\|_{(p^k)} \|V_{p^k} d^k + X(p^k) \|_{(p^k)}\\
        & \overset{\cref{eq: residual ineq}}{\leq} \|\log_{p^k} p^* + V_{p^k}^{-1} X(p^k)\|_{(p^k)} + a^{k}\|V_{p^k}^{-1}\|_{(p^k)} \|X(p^k)\|_{(p^k)}.
        \label{eq: split - thm inexact conv}
        \end{align}
        Since $X(p^*) = 0$ and parallel transport is an isometry, we see that
        \begin{align}
        \|\log_{p^k} p^* + V_{p^k}^{-1} X(p^k)\|_{(p^k)} & =  \|V_{p^k}^{-1} \left( V_{p^k} \log_{p^k} p^* +  X(p^k)\right)\|_{(p^k)} \\
        & \leq \|V_{p^k}^{-1} \|_{(p^k)} 
        \|\mathcal{P}_{p^{*} \leftarrow p^k}\left(X(p^k)+V_{p^k} \log _{p^k} p^{*} \right) \|_{(p^*)}\\
        & \leq \|V_{p^k}^{-1} \|_{(p^k)} 
        \|X\left(p^{*}\right)-\mathcal{P}_{p^{*} \leftarrow p^k} \left(X(p^k)+V_{p^k} \log _{p^k} p^{*} \right) \|_{(p^*)}.
        \label{eq: inexact step}
        \end{align}
        Substituting \cref{eq: inexact step} back into \cref{eq: split - thm inexact conv}, we find 
        \begin{align}
        \| \log_{p^k} p^* - d^k &\|_{(p^k)}\nonumber  \\
        \leq &  \|V_{p^k}^{-1} \|_{(p^k)}  \left(\|X\left(p^{*}\right)-\mathcal{P}_{p^{*} \leftarrow p^k} \left(X(p^k)+V_{p^k} \log _{p^k} p^{*} \right) \|_{(p^*)} +       a^k \|X(p^k)\|_{(p^k)}
        \right).
        \end{align}
        With the bounds derived earlier, we estimate
        \begin{align}
        \| \log_{p^k} p^* - d^k \|_{(p^k)} \overset{\cref{eq: V bound}, \cref{eq: semismooth}, \cref{eq: thm - lipschitz}}{\leq}&  \frac{ \lambda_{p^{*}}}{1-\epsilon \lambda_{p^{*}}} (\epsilon d_{\mathcal{M}}(p^k,p^*) + a^k L d_{\mathcal{M}}(p^k,p^*))\\
        \overset{a^k \leq a}{\leq}\hspace{1.75em} & \frac{ \lambda_{p^{*}}}{1-\epsilon \lambda_{p^{*}}} (\epsilon + a L) d_{\mathcal{M}}(p^k,p^*).
        \label{eq: vec leng res}
        \end{align}
        Now note that, since $K_{p^*} \geq 1$ (see \cref{rem: Kp}), we have 
        \begin{equation}
        \label{eq: thm K ineq}
        \frac{ \lambda_{p^{*}}}{1-\epsilon \lambda_{p^{*}}} (\epsilon + a L)  \leq \frac{ \lambda_{p^{*}}K_{p^*}}{1-\epsilon \lambda_{p^{*}}} (\epsilon + a L) .
        \end{equation}
        For our choice of $a$ and $\epsilon,$  we find
        \begin{align}
                &\epsilon \lambda_{p^{*}}\left(1+K_{p^{*}}\right)<1 - a \lambda_{p^{*}} L K_{p^{*}}\\
                \Leftrightarrow \quad & \epsilon \lambda_{p^{*}} + \epsilon \lambda_{p^{*}}K_{p^{*}} +  a \lambda_{p^{*}} L K_{p^{*}}<1 \\
                 \Leftrightarrow \quad &  \lambda_{p^{*}}K_{p^{*}}(\epsilon+ aL) < 1- \epsilon \lambda_{p^{*}} \\
                 \Leftrightarrow \quad & \frac{ \lambda_{p^{*}}K_{p^*}}{1-\epsilon \lambda_{p^{*}}} (\epsilon + a L)  < 1.
                 \label{eq: thm klp ineq}
        \end{align}
        Since $d_{\mathcal{M}}(p^k,p^*)<\delta$, we obtain from combining \cref{eq: vec leng res}, \cref{eq: thm K ineq} and \cref{eq: thm klp ineq}
        \begin{equation}
                \|\log_{p^k} p^* + V_{p^k}^{-1} X(p^k)\|_{(p^k)} < d_{\mathcal{M}}(p^k,p^*) <\delta \leq r_{p^*}.
        \end{equation}
         Moreover, we have $\|\log_{p^k}p^*\|_{(p^k)} = d_{\mathcal{M}}(p^k,p^*) \leq r_{p^*}$. Hence, we find (see \cref{def: Kp})
        \begin{equation}
        d_{\mathcal{M}}(\exp_{p^k}(d^k), p^*) \leq K_{p^*} \|\log_{p^k} p^* - d^k\|_{(p^k)}
        \label{eq: d K ineq}
        \end{equation}
    and observe
        \begin{align}
        \frac{d_{\mathcal{M}}(p^{k+1}, p^*)}{d_{\mathcal{M}}(p^k,p^*)} &\overset{\text{\cref{alg: inexact RSSN}}}{=} \frac{d_{\mathcal{M}}(\exp_{p^k}(d^k), p^*)}{d_{\mathcal{M}}(p^k,p^*)} \overset{\cref{eq: d K ineq}}{\leq} \frac{K_{p^*} \|\log_{p^k} p^* - d^k\|_{(p^k)}}{d_{\mathcal{M}}(p^k,p^*)} \\
        &\overset{\cref{eq: vec leng res}}{\leq} \frac{ \lambda_{p^{*}} K_{p^*}}{1-\epsilon \lambda_{p^{*}}} (\epsilon + a L) \overset{\cref{eq: thm klp ineq}}{<}1.
        \end{align}
        From this result we conclude by induction that if we choose $p^0 \in B_{\delta}(p^*)$ as in the assumption, we have $p^k \in B_{\delta}(p^*)$ $\forall k \in \mathbb{N}$ and convergence is Q-linear.

        (ii) The second part is very similar. Let $K_{p^*}$, $r_{p^*,}$ and  $\hat{\delta}$ with corresponding $L$ as before. Choose $\epsilon >0 $ such that $\epsilon \lambda_{p^*} (1 + 2K_{p^*})<1$ and take $0< \delta < \min\{\hat{\delta},r_{p^*}\} $ such that \cref{eq: V bound} and \cref{eq: semismooth} hold. Due to the assumption $\|r^k\|_{(p^k)}\in o\left(\|X(p^k)\|_{(p^k)}\right)$, the assumption that $p^k \rightarrow p^*$, and that $X$ is continuous, we have that $\|X(p^k)\|_{(p^k)} \rightarrow 0$. Moreover for large enough $k$,
        \begin{equation}
        \|r^k \|_{(p^k)} < \epsilon \delta.
        \label{eq: est r small}
        \end{equation}
        Because of the convergence assumption $p^k \rightarrow p^*$, we also have for large $k$ that $d_{\mathcal{M}}(p^k,p^*)<\delta$. Consequently, using the similar steps that lead to \cref{eq: vec leng res} in the proof of (i), we see
        \begin{align}
        \label{eq: ii ed r ineq}
        \| \log_{p^k} p^* - d^k \|_{(p^k)} & \leq \frac{ \lambda_{p^{*}}}{1-\epsilon \lambda_{p^{*}}} (\epsilon d_{\mathcal{M}}(p^k,p^*) + \|r^k \|_{(p^k)}) \\
        & \leq \frac{ 2 \epsilon\lambda_{p^{*}}}{1-\epsilon \lambda_{p^{*}}} \delta\leq \frac{ 2 \epsilon\lambda_{p^{*}}K_{p^{*}}}{1-\epsilon \lambda_{p^{*}}} \delta .
        \label{eq: ii eq log}
        \end{align}
        For our choice of $\epsilon$ we find
        \begin{align}
                &\epsilon \lambda_{p^{*}}\left(1+2K_{p^{*}}\right)<1 \\
                \Leftrightarrow \quad &  2\epsilon\lambda_{p^{*}}K_{p^{*}} < 1- \epsilon \lambda_{p^{*}} \\
                \Leftrightarrow \quad & \frac{2\epsilon \lambda_{p^{*}}K_{p^*}}{1-\epsilon \lambda_{p^{*}}}  < 1.
                \label{eq: ii ekl <1}
        \end{align}
        Since $d_{\mathcal{M}}(p^k,p^*)<\delta$, we obtain from combining \cref{eq: ii eq log} and \cref{eq: ii ekl <1}
        \begin{equation}
        \|\log_{p^k} p^* + V_{p^k}^{-1} X(p^k)\|_{(p^k)} <\delta \leq r_{p^*}.
        \end{equation}
        Again, we have $\|\log_{p^k}p^*\|_{(p^k)} = d_{\mathcal{M}}(p^k,p^*) \leq r_{p^*}$, which leads to (see \cref{def: Kp})
        \begin{equation}
        d_{\mathcal{M}}(\exp_{p^k}(d^k), p^*) \leq K_{p^*} \|\log_{p^k} p^* - d^k\|_{(p^k)}.
        \label{eq: de log d}
        \end{equation}
        Finally we see that, for large $k$,
        \begin{align}
        \frac{d_{\mathcal{M}}(p^{k+1}, p^*)}{d_{\mathcal{M}}(p^k,p^*)} &\overset{\text{\cref{alg: inexact RSSN}}}{=}\frac{d_{\mathcal{M}}(\exp_{p^k}(d^k), p^*)}{d_{\mathcal{M}}(p^k,p^*)} \overset{\cref{eq: de log d}}{\leq}\frac{K_{p^*} \|\log_{p^k} p^* - d^k\|_{(p^k)}}{d_{\mathcal{M}}(p^k,p^*)}\\
        &\overset{\cref{eq: ii ed r ineq}}{\leq} \frac{\epsilon\lambda_{p^{*}}K_{p^*}}{1-\epsilon \lambda_{p^{*}}} + \frac{\lambda_{p^{*}}K_{p^*}}{1-\epsilon \lambda_{p^{*}}} \frac{\|r^k\|_{(p^k)}}{d_{\mathcal{M}}(p^k,p^*)}.\\
        \intertext{By $\|X(p^k)\|_{(p^k)} \leq L d_{\mathcal{M}}(p^k,p^*)$, we can continue}
        &\leq\frac{\epsilon\lambda_{p^{*}}K_{p^*}}{1-\epsilon \lambda_{p^{*}}} + \frac{1}{L}\frac{\lambda_{p^{*}}K_{p^*}}{1-\epsilon \lambda_{p^{*}}} \frac{\|r^k\|_{(p^k)}}{\|X(p^k)\|_{(p^k)}}.
        \label{eq: ii el + rx res}
        \end{align}
        Note that this result holds for all (arbitrarily small) $\epsilon >0 $ such that $\epsilon \lambda_{p^*} (1 + 2K_{p^*})<1$. Hence, we can focus solely on the residual term. Then, by the assumption $\|r^k\|_{(p^k)}\in o(\|X(p^k)\|_{(p^k)})$,
        \begin{equation}
        \lim_{k \rightarrow \infty} \frac{d_{\mathcal{M}}(p^{k+1}, p^*)}{d_{\mathcal{M}}(p^k,p^*)} \leq \lim_{k \rightarrow \infty}\frac{\lambda_{p^{*}}K_{p^*}}{L(1-\epsilon \lambda_{p^{*}})} \frac{\|r^k\|_{(p^k)}}{\|X(p^k)\|_{(p^k)}} = 0,
        \end{equation}
        and we conclude that the convergence is superlinear.
        
        (iii) Again this is very similar to the previous cases. Let $\mu$ denote the order of the semi-smoothness and let $K_{p^*}$, $r_{p^*}$, and $\hat{\delta}$ with corresponding $L$ as before. For $\epsilon>0$ such that $ \epsilon \lambda_{p^*} <1$, choose $0< \delta < \min\{\hat{\delta},r_{p^*}\} $ satisfying $\epsilon\lambda_{p^*}(1 + 2\delta^\mu K_{p^*}) < 1$ and such that \cref{eq: V bound} and
        \begin{equation}
        \left\|X(p^*)-\mathcal{P}_{p^* \leftarrow p}\left[X(p)+V_{p} \exp _{p}^{-1} p^*\right]\right\|_{(p^*)} \leq \epsilon d_{\mathcal{M}}(p, p^*)^{1+\mu}
        \end{equation}
        hold. Then for large enough $k$ we have by the same reasoning as for establishing \cref{eq: est r small}, that
        \begin{equation}
        \|r^k \|_{(p^k)} < \epsilon \delta^{1+\mu}.
        \end{equation}
        Similar to \cref{eq: ii ekl <1},
        \begin{align}
        \label{eq: iii ed r ineq}
        \| \log_{p^k} p^* - d^k \|_{(p^k)} & \leq \frac{ \lambda_{p^{*}}}{1-\epsilon \lambda_{p^{*}}} (\epsilon d_{\mathcal{M}}(p^k,p^*)^{1+\mu} + \|r^k \|_{(p^k)}) \\
        & \leq \frac{ 2 \epsilon\lambda_{p^{*}}\delta^\mu}{1-\epsilon \lambda_{p^{*}}} \delta\leq \frac{ 2 \epsilon\lambda_{p^{*}}\delta^\mu K_{p^{*}}}{1-\epsilon \lambda_{p^{*}}} \delta
        \label{eq: iii eq log}
        \end{align}
        follows, and for our choice of $\epsilon$ and $\delta$ we find
        \begin{align}
        \epsilon \lambda_{p^{*}}\left(1+2\delta^\mu K_{p^{*}}\right)<1 
        \Leftrightarrow  \frac{2\epsilon \lambda_{p^{*}}\delta^\mu K_{p^*}}{1-\epsilon \lambda_{p^{*}}}  < 1.
        \label{eq: iii ekl <1}
        \end{align}
        Using similar arguments as for establishing \cref{eq: ii el + rx res} in (ii), we obtain
        \begin{align}
        \frac{d_{\mathcal{M}}(p^{k+1}, p^*)}{d_{\mathcal{M}}(p^k,p^*)} \leq\frac{\epsilon\lambda_{p^{*}}K_{p^*}}{1-\epsilon \lambda_{p^{*}}} + \frac{1}{L}\frac{\lambda_{p^{*}}K_{p^*}}{1-\epsilon \lambda_{p^{*}}} \frac{\|r^k\|_{(p^k)}}{\|X(p^k)\|_{(p^k)}^{1+\mu}}.
        \label{eq: iii el + rx res}
        \end{align}
        As $\epsilon$ can be arbitrarily small, we can again focus on the second term as in (ii). Finally, we see by the assumption $\|r^k\|_{(p^k)} = \mathcal{O}(\|X(p^k)\|_{(p^k)}^{1+\mu})$ that
        \begin{equation}
        \lim_{k \rightarrow \infty} \frac{d_{\mathcal{M}}(p^{k+1}, p^*)}{d_{\mathcal{M}}(p^k,p^*)^{1+\mu}} \leq \lim_{k \rightarrow \infty}\frac{\lambda_{p^{*}}K_{p^*}}{L(1-\epsilon \lambda_{p^{*}})} \frac{\|r^k\|_{(p^k)}}{\|X(p^k)\|_{(p^k)}^{1+\mu}} =: M
        \end{equation}
        for some $M>0$. Therefore we conclude that the convergence is of Q-order $1+\mu$.
\end{proof}

From (ii) we obtain the following practical condition if one can ensure that the relative error in solving the Newton system approaches zero.
\begin{corollary}
        If the sequence $(p^{k})_{k\geq 0}$ generated by \cref{alg: inexact RSSN} converges to the solution $p^{*}$ and the sequence $\{a^k\}$ converges to zero, then the rate of convergence is $Q$-superlinear.
\end{corollary}

\begin{remark}
        In the real-valued case (ii) and (iii) are formulated stronger: if a convergent sequence exists, the converse statements in (ii) and (iii) also hold. However, the proof in \cite{facchinei1996inexact} relies heavily on the linearity of $\mathbb{R}^d$ and is therefore much harder to translate to the manifold case. This remains an open problem.
\end{remark}

\subsection{Observations}
\label{sec: observations rssn}
Before moving on to applications, we will elaborate on some key observations in the theory of the (inexact) RSSN method. 
In particular, we can make the following observations regarding the convergence behavior of both algorithms. By a continuity argument, we see that the minimal radius of convergence $\delta$ in \cref{thm: inexact Newton} is determined by~$\epsilon$. This $\epsilon$ must satisfy the upper bound (that is for $a^k = 0$)
\begin{equation}
\epsilon \lambda_{p^{*}}\left(1+K_{p^{*}}\right)<1.
\end{equation}
\jla{We see that we run into trouble in two cases: large $\lambda_{p^*}$ and large $K_{p^{*}}$. This comes down to the following scenarios:}

\subparagraph{The generalized covariant derivative is close to singular}
Here a large $\lambda_{p^{*}}$ only admits a very small $\epsilon$, which in turn results in a small convergence region.

\subparagraph{The manifold has negative curvature}
For negatively-curved manifolds we do not have an \emph{a priori} estimate for $K_{p^*}$. If this value becomes arbitrarily large, we cannot expect a large region of convergence. 



\section{Application to $\ell^2$-TV-like Functionals}
\label{sec: application}


We consider the classical Rudin-Osher-Fatemi denoising problem, extended to the manifold-valued setting:
Let $\mathcal{M}$ be a Riemannian manifold, $d_1,d_2 \in \mathbb{N}$ be the dimensions of the image, and $h \in \mathcal{M}^{d_1\times d_2}$ be the noisy input data. We are interested in solving the isotropic ($q=2$) and anisotropic ($q = 1$) discrete ROF model \cite{bergmann2019fenchel}
\begin{equation}
\inf_{p \in \mathcal{M}^{d_1\times d_2}} \frac{1}{2 \alpha} \sum_{i, j=1}^{d_{1}, d_{2}} d_{\mathcal{M}}\left(p_{i, j},h_{i, j}\right)^{2} + \|T (p) \|_{p,q,1},
\label{eq: l2TV manifold}
\end{equation}
where $\alpha > 0 $ and $T: \mathcal{M}^{d_1\times d_2} \rightarrow \mathcal{T}\mathcal{M}^{d_1\times d_2 \times 2}$ is the non-linear finite difference operator~{\cite{bergmann2019fenchel}}
\begin{equation}
(T (p))_{i, j, k}:=\left\{\begin{array}{ll}
0 \in \mathcal{T}_{p_{i, j}} \mathcal{M} & \text { if } i=d_{1} \text { and } k=1 \\
0 \in \mathcal{T}_{p_{i, j}} \mathcal{M} & \text { if } j=d_{2} \text { and } k=2 \\
\log _{p_{i, j}} p_{i+1, j} \in \mathcal{T}_{p_{i, j}} \mathcal{M} & \text { if } i<d_{1} \text { and } k=1 \\
\log _{p_{i, j}} p_{i, j+1} \in \mathcal{T}_{p_{i, j}} \mathcal{M} & \text { if } j<d_{2} \text { and } k=2
\end{array}\right.
\end{equation}
and with the norm defined as
\begin{equation}
\|T (p) \|_{p,q,1} := \sum_{i,j = 1}^{d_1,d_2} (\|(T (p))_{i,j,1}\|_{p_{i,j}}^q + \|(T (p))_{i,j,2}\|_{p_{i,j}}^q)^{\frac{1}{q}}.
\end{equation}
The first step is to apply Fenchel duality theory in order to rewrite  model \cref{eq: l2TV manifold} into an appropriate linearized saddle-point form with corresponding optimality conditions. We introduce the following notation:

First, as $T(p) \in \mathcal{T}_p\mathcal{M}^{d_1\times d_2}\times \mathcal{T}_p\mathcal{M}^{d_1\times d_2} \cong \mathcal{T}_{(p,p)}\mathcal{M}^{d_1\times d_2\times 2}$, we write the base point $(p,p)$ of $T(p)$ simply as $p$,  and consequently write $\mathcal{T}_{p}\mathcal{M}^{d_1\times d_2\times 2}$ for $T(p)$ as well.

Next, the dual space $\mathcal{T}^*_n \mathcal{T}\mathcal{M}^{d_1\times d_2\times 2}$ for base point $n \in \mathcal{T}\mathcal{M}^{d_1\times d_2\times 2}$ can be constructed through the \emph{double tangent bundle} $\mathcal{T}^2\mathcal{M}^{d_1\times d_2\times 2} := \mathcal{T}\mathcal{T}\mathcal{M}^{d_1\times d_2\times 2} $. Elements in the latter can be written as $(p,X_p,Y_p,Z_p) \in \mathcal{T}^2\mathcal{M}^{d_1\times d_2\times 2}$, where $p\in \mathcal{M}^{d_1\times d_2\times 2}$ and $X_p,Y_p,Z_p \in \mathcal{T}_p\mathcal{M}^{d_1\times d_2\times 2}$ {\cite[Sec. 2]{li2013geometry}}. We also write $(Y_p,Z_p ) \in \mathcal{T}_{(p,X_p)} \mathcal{T} \mathcal{M}^{d_1\times d_2\times 2}$. Similarly, for elements in the dual bundle we write $(p,X_p,\eta_p,\xi_p)  \in \mathcal{T}^* \mathcal{T}\mathcal{M}^{d_1\times d_2\times 2}$ and $(\eta_p,\xi_p ) \in \mathcal{T}^*_{(p,X_p)} \mathcal{T} \mathcal{M}^{d_1\times d_2\times 2} $, where $p\in \mathcal{M}^{d_1\times d_2\times 2}$, $X_p \in \mathcal{T}_p\mathcal{M}^{d_1\times d_2\times 2}$ and $\eta_p,\xi_p  \in \mathcal{T}^*_p \mathcal{M}^{d_1\times d_2\times 2}$. For the duality pairing we have $\langle (\eta_p,\xi_p), (Y_p,Z_p ) \rangle_{(p,X_p)} := \langle \eta_p, Y_p \rangle_p + \langle \xi_p, Z_p \rangle_p$.

In the following, we will need the logarithmic mapping on the tangent bundle $\mathcal{T}\mathcal{M}^{d_1\times d_2\times 2}$. The exponential mapping  is given by $\exp_{(p,X_p)} (Y_p, Z_p ) := (\exp_p Y_p, \mathcal{P}_{\exp_p Y_p \leftarrow p} (Z_p + X_p)) \in \mathcal{T}\mathcal{M}^{d_1\times d_2\times 2}$ \cite[p. 30]{bergmann2019fenchel}. Hence, for the logarithmic mapping we must have $\log_{(p,X_p)} (q,Y_q) := (\log_p q, \mathcal{P}_{p \leftarrow q}Y_q- X_p) \in \mathcal{T}_{(p,X_p)} \mathcal{T}\mathcal{M}^{d_1\times d_2\times 2}$.

Then, using duality of the $\|\cdot \|_{p,q,1}$ norm and the fact that parallel transport is an isometry, we can choose $m \in \mathcal{M}$ and write 
\begin{equation}
        \|T (p) \|_{p,q,1} = \|\mathcal{P}_{m \leftarrow p} T (p) \|_{m,q,1} = \sup_{\xi_m \in T_{m}^*\mathcal{M}^{d_1\times d_2\times 2}} \langle \xi_{m}, \mathcal{P}_{m \leftarrow p} T (p) \rangle_m - \iota_{B_{q^*}}(\xi_m),
\end{equation}
where $q^*$ such that $\frac{1}{q} + \frac{1}{q^*} = 1$,
\begin{align}
B_{q^{*}} &:= \left\{\nu_m \in  \mathcal{T}_{m}^*\mathcal{M}^{d_1\times d_2\times 2} \mid\|\nu_m\|_{m,q^{*}, \infty} \leq 1\right\}\\
& = \left\{\nu_m \in   \mathcal{T}_{m}^*\mathcal{M}^{d_1\times d_2\times 2} \mid \max \|(\nu_m)_{i,j,:}\|_{m,q^{*}} \leq 1\right\} ,
\end{align}
and 
\begin{equation}
\iota_{B_{q^*}}(\xi_m) := \left\{\begin{array}{ll}
0  & \text { if } \xi_m \in B_{q^*}, \\
\infty  & \text { if } \xi_m \notin B_{q^*}.
\end{array}\right. 
\end{equation}
This choice of $m$ can be seen as a choice of origin on the manifold: while in classical duality theory on vector spaces we use the zero element, manifolds do not have such a preferred point. We will discuss the specific choice of $m$ for this application at the end of the section.

First, we choose $n := 0 \in \mathcal{T}_{m}\mathcal{M}^{d_1\times d_2\times 2}$ as the zero vector in the tangent space, and rewrite
\begin{align}
        \langle  \xi_{m} , \mathcal{P}_{m \leftarrow p} T (p)\rangle_m &= \langle \xi_{m} , \mathcal{P}_{m \leftarrow p} T (p) - n\rangle_m \\
        & = \langle  \xi_{m}, \mathcal{P}_{m \leftarrow p} T (p) - n\rangle_m  + \sup_{\eta_m \in T_{m}^*\mathcal{M}^{d_1\times d_2\times 2}} \{  \langle \eta_m , \log_m p\rangle_m - \iota_{\{0\}}(\eta_m) \}\\
        & = \sup_{\eta_m \in T_{m}^*\mathcal{M}^{d_1\times d_2\times 2}} \{ \langle (\eta_m,\xi_{m}) , (\log_m p, \mathcal{P}_{m \leftarrow p} T (p)-n)\rangle_n - \iota_{\{0\}}(\eta_m) \}\\
        & = \sup_{\eta_m \in T_{m}^*\mathcal{M}^{d_1\times d_2\times 2}} \{ \langle (\eta_m,\xi_{m}), \log_n T (p) \rangle_n - \iota_{\{0\}}(\eta_m)\},
\end{align}
where we added $\sup_{\eta_m \in T_{m}^*\mathcal{M}^{d_1\times d_2\times 2}}   \{\langle \eta_m , \log_m p\rangle_m - \iota_{\{0\}}(\eta_m)\}=0$ in the second line.

Finally, let $\xi_n := (\eta_m, \xi_m) \in \mathcal{T}_n^* \mathcal{T}\mathcal{M}^{d_1\times d_2 \times 2}$. Bringing everything together, we find a saddle-point problem in the desired form:
\begin{equation}
\inf_{p \in \mathcal{M}^{d_1\times d_2}} \sup_{\xi_{n} \in \mathcal{T}_{n}^* \mathcal{T}\mathcal{M}^{d_1\times d_2\times 2}}  \frac{1}{2 \alpha} \sum_{i, j=1}^{d_{1}, d_{2}} d_{\mathcal{M}}\left(p_{i, j},h_{i, j}\right)^{2} + \langle \xi_n, \log_n T(p)  \rangle_n - \iota_{\{0\}\times B_{q^*}}(\xi_n).
\end{equation}
Choosing $m$ such that $T(m)=n$, the corresponding linearized saddle-point problem is
\begin{equation}
        \inf_{p \in \mathcal{M}^{d_1\times d_2}} \sup_{\xi_{n} \in \mathcal{T}_{n}^* \mathcal{T}\mathcal{M}^{d_1\times d_2\times 2}}  \frac{1}{2 \alpha} \sum_{i, j=1}^{d_{1}, d_{2}} d_{\mathcal{M}}\left(p_{i, j},h_{i, j}\right)^{2} + \left\langle \xi_{n}, D_m T\left[\log _{m} p\right] \right\rangle_n - \iota_{\{0\}\times B_{q^*}}(\xi_n),
        \label{eq: TV lin saddle}
\end{equation}
on which we can now apply the PD-RSSN method from \cref{sec: higher-order}.
Finally, it remains to find values for $m$ such that both $T(m)=n$ and $n=0$ hold. 
We must have
\begin{equation}
        m_{ij} =\tilde{m} \qquad \text{ for } i = 1,\ldots,d_1\text{ and } j=1,\ldots, d_2      
\end{equation}
for some fixed $\tilde{m}\in \mathcal{M}$\jla{, i.e., $m$ is constant.}


\paragraph{Adaptations to the optimization problem}
There is one serious numerical issue with solving \cref{eq: TV lin saddle} using a PD-RSSN approach: The generalized covariant derivative can be non-invertible. For the real-valued case this is the case and can be solved by adding dual regularization {\cite[Sec. 3.4.2 \& 3.4.3]{diepeveen2020nonsmooth}}. We use a similar strategy in the manifold case: We add a quadratic penalty term for the duals $\xi_n$,
\begin{equation}
        \inf_{p \in \mathcal{M}^{d_1\times d_2}} \sup_{\xi_{n} \in \mathcal{T}_{n}^* \mathcal{T}\mathcal{M}^{d_1\times d_2\times 2}}  \frac{1}{2 \alpha} \sum_{i, j=1}^{d_{1}, d_{2}} d_{\mathcal{M}}\left(p_{i, j},h_{i, j}\right)^{2} + \left\langle \xi_{n}, D_m T\left[\log _{m} p\right] \right\rangle_n - \iota_{\{0\}\times B_{q^*}}(\xi_n) - \frac{\beta}{2}\|\xi_n\|_n^2,
        \label{eq: tv reg final saddle}
\end{equation}
where $0<\beta \ll 1$ to ensure that the solution to \cref{eq: tv reg final saddle} is as close as possible to that of \cref{eq: TV lin saddle}. Next, in order to make the numerics easier and more efficient, we use that
the first entry of $\xi_n$ will be zero to get a smaller matrix representation of the generalized covariant derivative, i.e., we solve
\begin{equation}
        \inf_{p \in \mathcal{M}^{d_1\times d_2}} \sup_{\xi_{m} \in \mathcal{T}^*_{m}\mathcal{M}^{d_1\times d_2\times 2}}  \frac{1}{2 \alpha} \sum_{i, j=1}^{d_{1}, d_{2}} d_{\mathcal{M}}\left(p_{i, j},h_{i, j}\right)^{2} + \left\langle \xi_{m}, \nabla_{(\cdot)} T \left[\log _{m} p\right] \right\rangle_m - \iota_{B_{q^*}}(\xi_m) - \frac{\beta}{2}\|\xi_m\|_m^2.
\end{equation}
Here we use the fact that the value of the \wdpa{second entry} of the differential of $T$ is that of the covariant derivative. We obtain the reduced optimality conditions
\begin{align}
        p &= \operatorname{prox}_{\sigma F}\left(\exp _{p}\left(\mathcal{P}_{p \leftarrow m} \left(-\sigma \left(\nabla_{(\cdot)} T \right)^{*}\left[ \xi_{m}\right]\right)^{\sharp}\right)\right), \\
        \xi_{m} &= \operatorname{prox}_{\tau G_{m,q}^{*}}\left(\xi_{m}+\tau\left(\nabla_{(\cdot)} T \left[\log _{m} p\right]\right)^{\flat}\right),
\end{align}
and the reduced vector field $X: \mathcal{M}^{d_1\times d_2} \times  \mathcal{T}_m^* \mathcal{M}^{d_1\times d_2 \times 2} \rightarrow \mathcal{T}_p\mathcal{M}^{d_1\times d_2} \times \mathcal{T}_m^* \mathcal{M}^{d_1\times d_2 \times 2}$
\begin{equation}
        X(p,\xi_m) := \begin{pmatrix}
                -\log_p \operatorname{prox}_{\sigma F}\left(\exp _{p}\left(\mathcal{P}_{p \leftarrow m} \left(-\sigma \left(\nabla_{(\cdot)} T \right)^{*}\left[ \xi_{m}\right]\right)^{\sharp}\right)\right)\\ 
                \xi_{m} - \operatorname{prox}_{\tau G_{m,q}^{*}}\left(\xi_{n}+\tau\left(\nabla_{(\cdot)} T \left[\log _{m} p\right]\right)^{\flat}\right)
        \end{pmatrix},
\end{equation}
where
\begin{equation}
F(p) := \sum_{i, j=1}^{d_{1}, d_{2}} d_{\mathcal{M}}\left(p_{i, j},h_{i, j}\right)^{2} \quad \text{ and } \quad G^*_{m,q} (\xi_m) := \iota_{B_{q^*}}(\xi_m) + \frac{\beta}{2}\|\xi_m\|_m^2.
\end{equation}
For the expressions of the proximal maps and the differentials, we again refer the reader to the thesis \cite[Sec. 6.5]{diepeveen2020nonsmooth}. The differential of $T$ and its adjoint can be found in \cite[Sec. 5]{bergmann2019fenchel}.

\jla{Unlike the real-valued case mentioned above, the uniform invertibility of the Newton operator for the regularized $\ell^2$-TV problem on arbitrary manifolds is still an open question, and we have to point to the numerical evidence outlined in the following section.}

In order to show semismoothness, we only sketch the proof and refer to \cite{diepeveen2020nonsmooth} for details. We already have that $\operatorname{prox}_{\sigma F}(x)$ is semismooth since it is smooth: by the generalized Taylor series in \cite{dedieu2003newton} we find semismoothness according to \cref{def: semismooth manifold}. The semismoothness of $\operatorname{prox}_{\tau G^*_q}(y)$ follows from invoking \cite[Prop. 7.4.6]{facchinei2007finite} as in the non-manifold case. Indeed we can do this, since the dual variable just lives in a finite dimensional vector space, which is isomorphic with $\mathbb{R}^d$. Hence we are justified to use PD-RSSN for $\ell^2$-TV. 


\section{Numerical Experiments}
\label{sec: numerics}



In this section we will explore the behavior of the Riemannian Semi-smooth Newton method in several numerical experiments with the ROF\ model on both \wdpa{the 2-sphere $S^2$ and the manifold of symmetric positive definite $3 \times 3$ matrices $\mathcal{P}(3)$}. This choice of manifolds allows to observe the behavior on both positively as well as negatively curved manifolds.

The proposed algorithms were implemented in Julia version 1.3.0 using \manopt \cite{Bergmann2019Manopt} and evaluated on a 2.4 GHz Intel Core i7 with 8 GB RAM. The source code is available at
\vspace{.5em}
\begin{center}
\href{https://github.com/wdiepeveen/Primal-Dual-RSSN}{\texttt{https://github.com/wdiepeveen/Primal-Dual-RSSN}}.
\end{center}
\vspace{.5em}
In order to verify basic performance of the complete method, we first consider a 1D problem with known minimizer. Then we provide a detailed runtime analysis and comparison to lRCPA~\cite{bergmann2019fenchel}. Finally, we numerically validate the convergence rates predicted by \cref{thm: inexact Newton}.


Throughout this section we will use the relative (residual) error
\begin{equation}
\epsilon^k_{rel} := \frac{\|X(p^k,\xi_m^k)\|_{(p^k,\xi_m^k)}}{\|X(p^0,\xi_m^0)\|_{(p^0,\xi_m^0)}}
\end{equation}
as measure for convergence. \wdpa{All time measurements provided are CPU times, measured using the CPUtic and CPUtoc commands of the \texttt{CPUTime.jl} package. }

\subsection{Signal with Known Minimizers}
In this first experiment we investigate the progression of the relative error as well as the distance to the exact solution for a problem with known minimizer. We 
consider the 1-dimensional piecewise constant signal from~\cite{bergmann2019fenchel}:
\begin{equation}
\label{eq: piecewise signal manifold}
h \in \mathcal{M}^{2\ell} \quad h_{i}:=\left\{\begin{array}{ll}
\hat{p}_1 & \text { if } i \leq \ell \\
\hat{p}_2 & \text { if } i>\ell
\end{array}\right. .
\end{equation}
For the $\ell^2$-TV problem as in \cref{eq: l2TV manifold} with this signal, the exact minimizer is known: for $\alpha > 0$ and $\hat{p}_1,\hat{p}_2 \in \mathcal{M}$, the minimizer $p^*$ is given by
\begin{equation}
p^*_{i}:=\left\{\begin{array}{ll}
p^*_1  & \text { if } i \leq \ell \\
p^*_2 & \text { if } i> \ell
\end{array}\right. ,
\end{equation}
where
\begin{equation}
\label{eq: exact sol}
        p^*_{1}=\gamma_{\hat{p}_{1}, \hat{p}_{2}}(\delta), \quad \quad  p^*_{2}=\gamma_{\hat{p}_{2}, \hat{p}_{2}}(\delta), \quad \text { and } \quad \delta=\min \left\{ \frac{1}{2},\frac{\alpha}{\ell}\frac{1}{d_{\mathcal{M}}\left(\hat{p}_{1}, \hat{p}_{2}\right)}\right\}.
\end{equation}
A proof can be found in \cite[Appendix 2]{diepeveen2020nonsmooth}

Furthermore, since $d_2 = 1$, the operator $T$ reduces to the 1-dimensional non-linear difference operator, i.e., we have $T:\mathcal{M}^{2\ell} \rightarrow \mathcal{T}\mathcal{M}^{2 \ell} $. This also means that the isotropic ($q=2$) and anisotropic ($q = 1$) cases reduce to the same functional.


For our tests we set $\beta=0$, $\ell = 10$, $\alpha = 5,$ and $\sigma = \tau = \frac{1}{2}$. For the stopping criterion, we used $\epsilon_{rel} = 10^{-10}$. Furthermore, for a given primal base point $m \in \mathcal{M}$ we chose $n = T(m) \in \mathcal{T}\mathcal{M}^{2\ell}$ for the dual base point. In order to satisfy the condition that $n$ is the zero vector in $\mathcal{T}_m\mathcal{M}^{2\ell}$, it is enough to choose $m_i$ the same in every grid point. 

For both manifolds we investigated the influence of the starting point. For the \enquote{cold start} we set $p^0:=h$, $\xi^0_n := 0$. For the \enquote{warm start} we again set $p^0:=h$, but modified the dual start to minor change $\xi_n^0: = \operatorname{prox}_{\tau G_{m,1}^{*}}\left(\tau\left(\nabla_{(\cdot)} T \left[\log _{m} p^0\right]\right)^{\flat}\right)$, the result of a single lRCPA dual step.

In the following we denote by $\tilde{p}$ the solution obtained from our method. The distance of the obtained solution $\tilde{p}$ to the exact solution $p^\ast$ are summarized in \cref{tab: distances}. 

\begin{table}[h!]
        \centering
        \begin{tabular}{ccccccc}
                \toprule
                & \multicolumn{3}{c}{Cold start} & \multicolumn{3}{c}{Warm start} \\ 
                \cmidrule(l){2-4} \cmidrule(l){5-7} 
                $\mathcal{M}$ &  $d_{\mathcal{M}}(p^*,\tilde{p})$ & \# Iterations & Time (s)& $d_{\mathcal{M}}(p^*,\tilde{p})$ & \# Iterations & Time (s)\\ 
                \midrule
                $S^2$ & 3.295 & $\geq$50 & 3.11 &  \textbf{0.0}  & 2 & 0.469  \\ 
                $\mathcal{P}(3)$ & 2.129e-13 & 1&  4.61 &  \textbf{1.838e-13} & 1 & 4.202 \\ \bottomrule
        \end{tabular}
        \caption{Distance between the exact minimizer $p^\ast$  of the 1D piecewise constant $\ell^2$-TV problem and the final iterate $\tilde{p}$ of our method on the $S^2$ and the $\mathcal{P}(3)$ manifold. Using a warm start with a rough estimate for primal and dual variables, the result is very close to the exact  minimizer. \wdpa{For the $S^2$ problem with cold start, PD-RSSN was terminated after 50 iterations.}}
        \label{tab: distances}
\end{table}

\paragraph{Case 1: $\mathcal{M} = S^2$}
We chose
\begin{equation}
        \hat{p}_{1}:=\frac{1}{\sqrt{2}}(1,1,0)^{\top}, \quad \hat{p}_{2}:=\frac{1}{\sqrt{2}}(1,-1,0)^{\top}
\end{equation}
and $m_i := (1,0,0)^{\top}$ for all $ i = 1,...,2\ell$.


With a cold start, PD-RSSN \wdpa{converges to a fixed point of the Newton iteration} and is terminated after 50 iterations in 3.11 seconds. \wdpa{That is, the Newton step gives a non-zero direction $d^k \in \mathcal{T}_p(S^2)^{2\ell} \times \mathcal{T}_m^* (S^2)^{2\ell}$ such that $p^{k+1} = p^k$ and $\xi_n^{k+1} = \xi_n^k$. This is possible for $d^k = (z^k,0) \in \mathcal{T}_p(S^2)^{2\ell} \times \mathcal{T}_m^* (S^2)^{2\ell}$, where $z^k \in \mathcal{T}_p(S^2)^{2\ell}$ consists of tangent vectors in $\mathcal{T}_{p_i} S^2$ with length equal to some multiple of $2\pi$. The iterates indicate that this is indeed what happens in the cold start case}. With the warm start, convergence to the solution of the linearized system is achieved after only two iterations in 0.469 seconds total. As it turns out this is the same point as the ROF minimizer, i.e., $d_{\mathcal{M}}(p^*,\tilde{p})=0$ as well.

 \begin{figure}[h!]
        \centering
        \begin{subfigure}{0.49\textwidth}
                \centering
                \includegraphics[width=0.9\linewidth]{"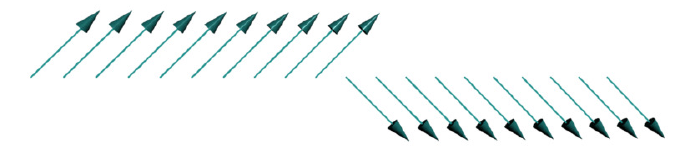"}
                \caption{Input}
        \end{subfigure}
        \hfill
        \begin{subfigure}{0.49\textwidth}
                \centering
                 \includegraphics[width=0.9\linewidth]{"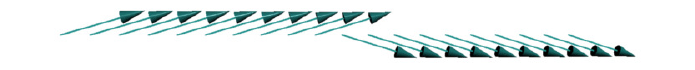"}
                \caption{Exact minimizer}
        \end{subfigure}\\
        \begin{subfigure}{0.49\textwidth}
                \centering
                \includegraphics[width=0.9\linewidth]{"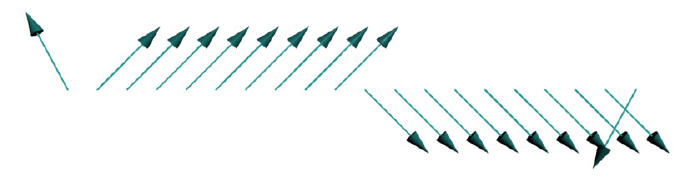"}
                \caption{Cold-start PD-RSSN solution}
        \end{subfigure}
        \hfill
        \begin{subfigure}{0.49\textwidth}
        \centering
        \includegraphics[width=0.9\linewidth]{"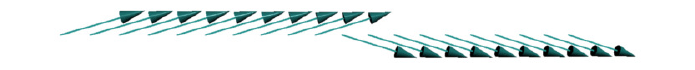"}
        \caption{Warm-start PD-RSSN solution}
 \end{subfigure}\\
        \begin{subfigure}{0.49\textwidth}
                \centering
                \includegraphics[width=0.99\linewidth]{"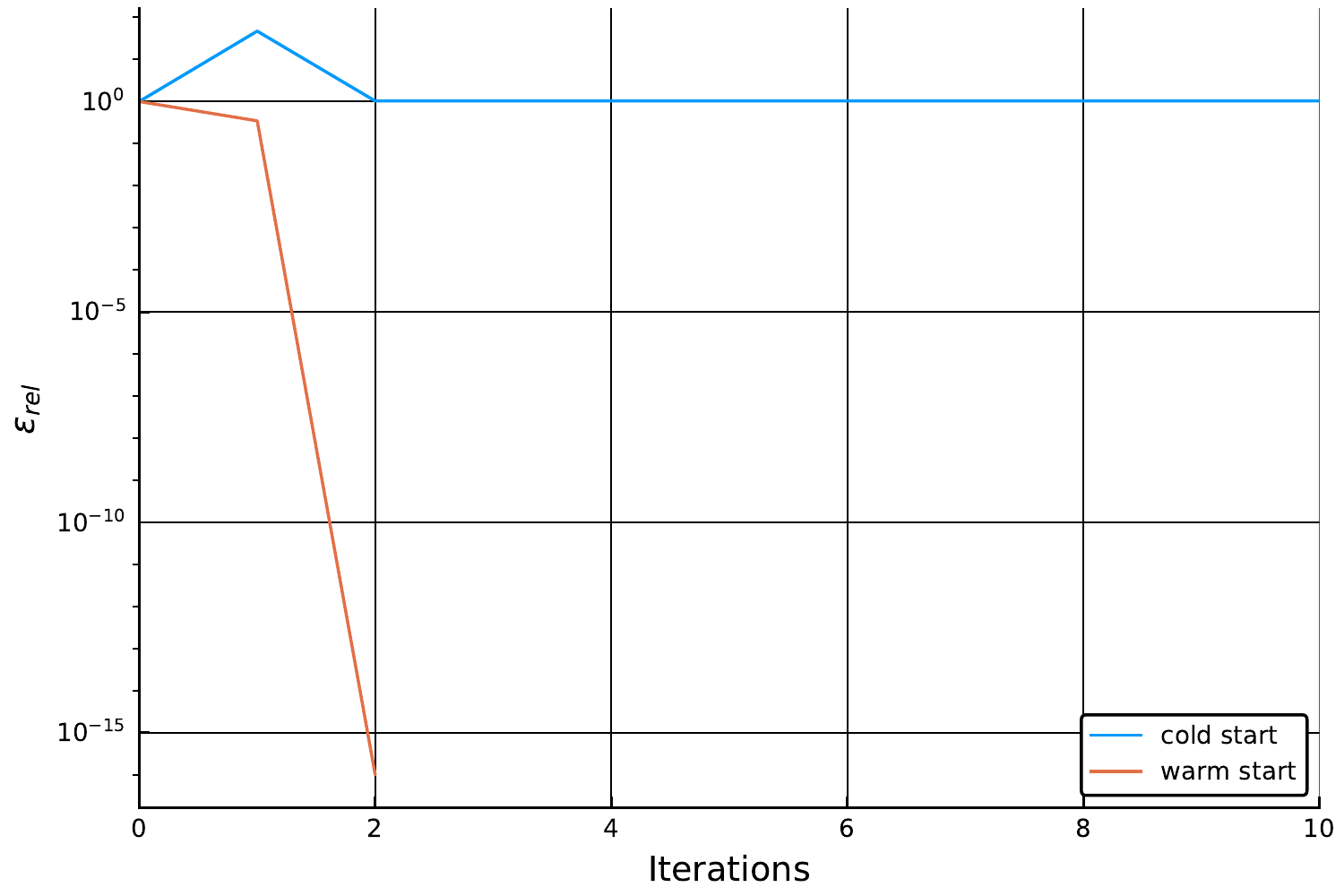"}
                \caption{Progression of the relative error}
        \end{subfigure}
        \hfill
        \begin{subfigure}{0.49\textwidth}
                \centering
                \includegraphics[width=0.99\linewidth]{"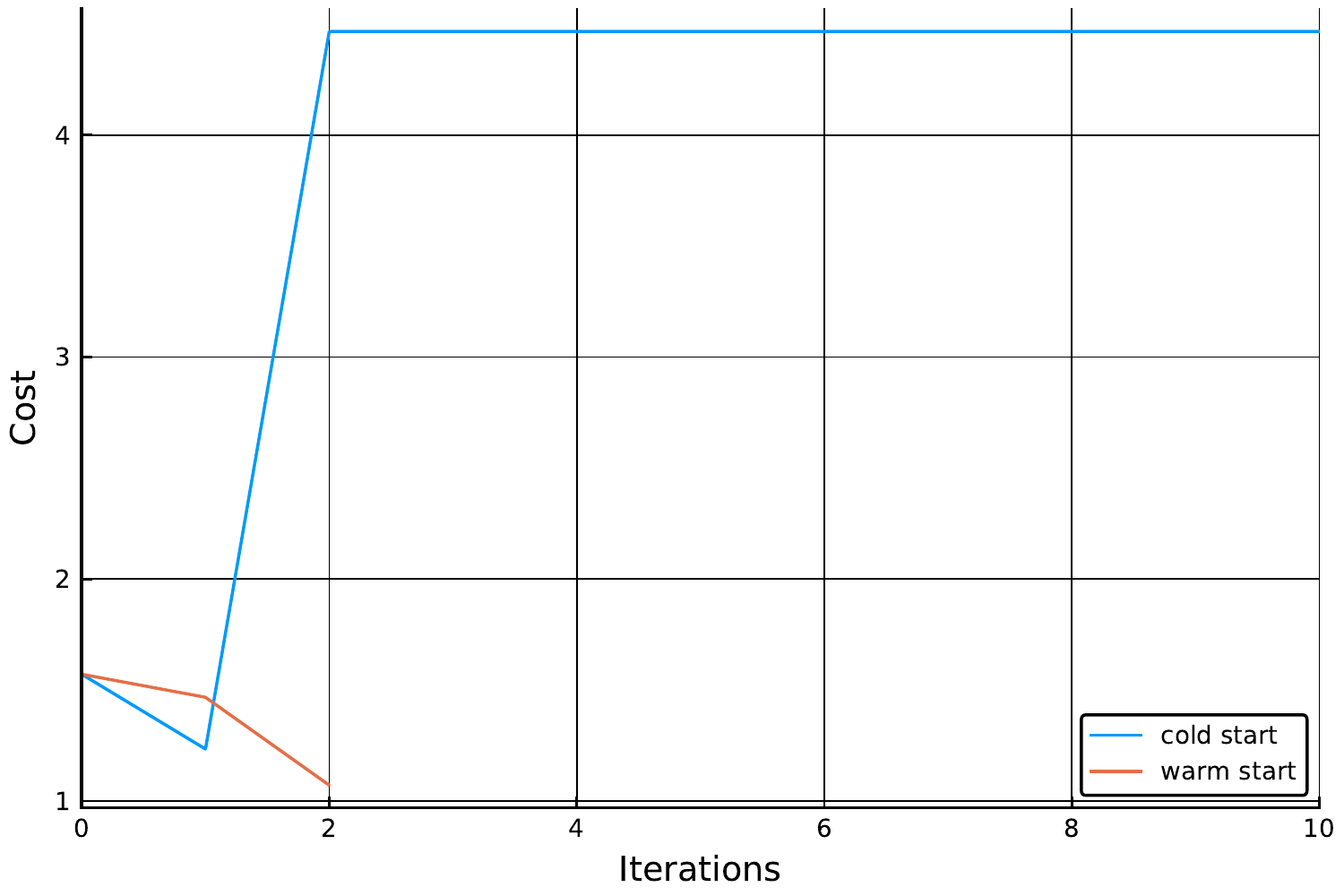"}
                \caption{Progression of the $\ell^2$-TV cost}
        \end{subfigure}\\
        \caption{
                Results and progression of the relative error and $\ell^2$-TV cost for the proposed PD-RSSN method for solving an  $S^2$-valued problem with known minimizer. With a cold start, PD-RSSN \wdpa{converges to a fixed point of the Newton iteration}. With a warm-start strategy, the algorithm converges rapidly to the exact minimizer.
        }
        \label{fig: RSSN-known-minimizers-S2}
 \end{figure}

\paragraph{Case 2: $\mathcal{M} = \mathcal{P}(3)$}
Here we chose
\begin{equation}
        \hat{p}_{1}:=\exp _{I}\left(\frac{2}{\|X\|_{I}} X\right), \quad \hat{p}_{2}:=\exp _{I}\left(-\frac{2}{\|X\|_{I}} X\right),\quad  \text { with } \quad X:=\left(\begin{array}{ccc}
        1 & 2 & 2 \\
        2 & 2 & 0 \\
        2 & 0 & 6
        \end{array}\right) ,
\end{equation} 
where $X \in \mathcal{T}_{I} \mathcal{P}(3)$ and $I$ is the identity matrix, and pick $m_i := I $ for all $i = 1,...,2\ell$.


In this setting,  we observed no difference between warm start and cold start (\cref{tab: distances}). With the cold start, PD-RSSN converges in 1 iteration in 4.61 seconds at a distance of $d_{\mathcal{M}}(p^*,p_l) = 2.13\cdot 10^{-13}$ from the exact solution. With the warm start, PD-RSSN converges after 1 iteration in 4.202 seconds at an almost identical distance of $d_{\mathcal{M}}(p^*,p_l) =1.84\cdot 10^{-13}$.

 \begin{figure}[h!]
        \centering
        \begin{subfigure}{0.49\textwidth}
                \centering
                \includegraphics[width=0.9\linewidth]{"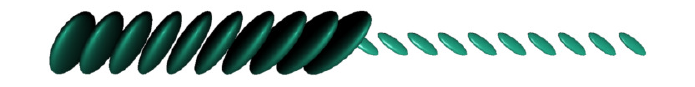"}
                \caption{Input}
        \end{subfigure}
        \hfill
        \begin{subfigure}{0.49\textwidth}
                \centering
                 \includegraphics[width=0.9\linewidth]{"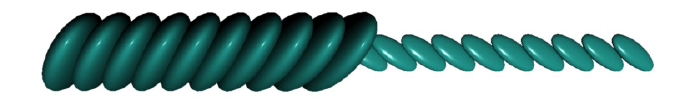"}   
                \caption{Exact minimizer}
        \end{subfigure}\\
        \begin{subfigure}{0.49\textwidth}
                \centering
                 \includegraphics[width=0.9\linewidth]{"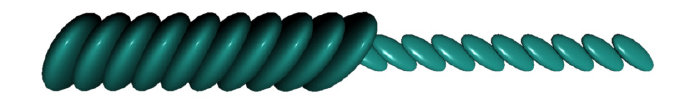"}
                \caption{Cold-start PD-RSSN solution}
        \end{subfigure}
        \hfill
        \begin{subfigure}{0.49\textwidth}
        \centering
        \includegraphics[width=0.9\linewidth]{"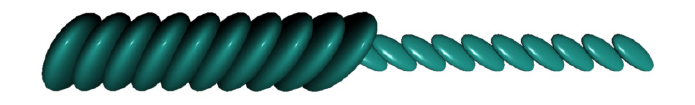"}
        \caption{Warm-start PD-RSSN solution}
\end{subfigure} \\
        \begin{subfigure}{0.49\textwidth}
                \centering
                \includegraphics[width=\linewidth]{"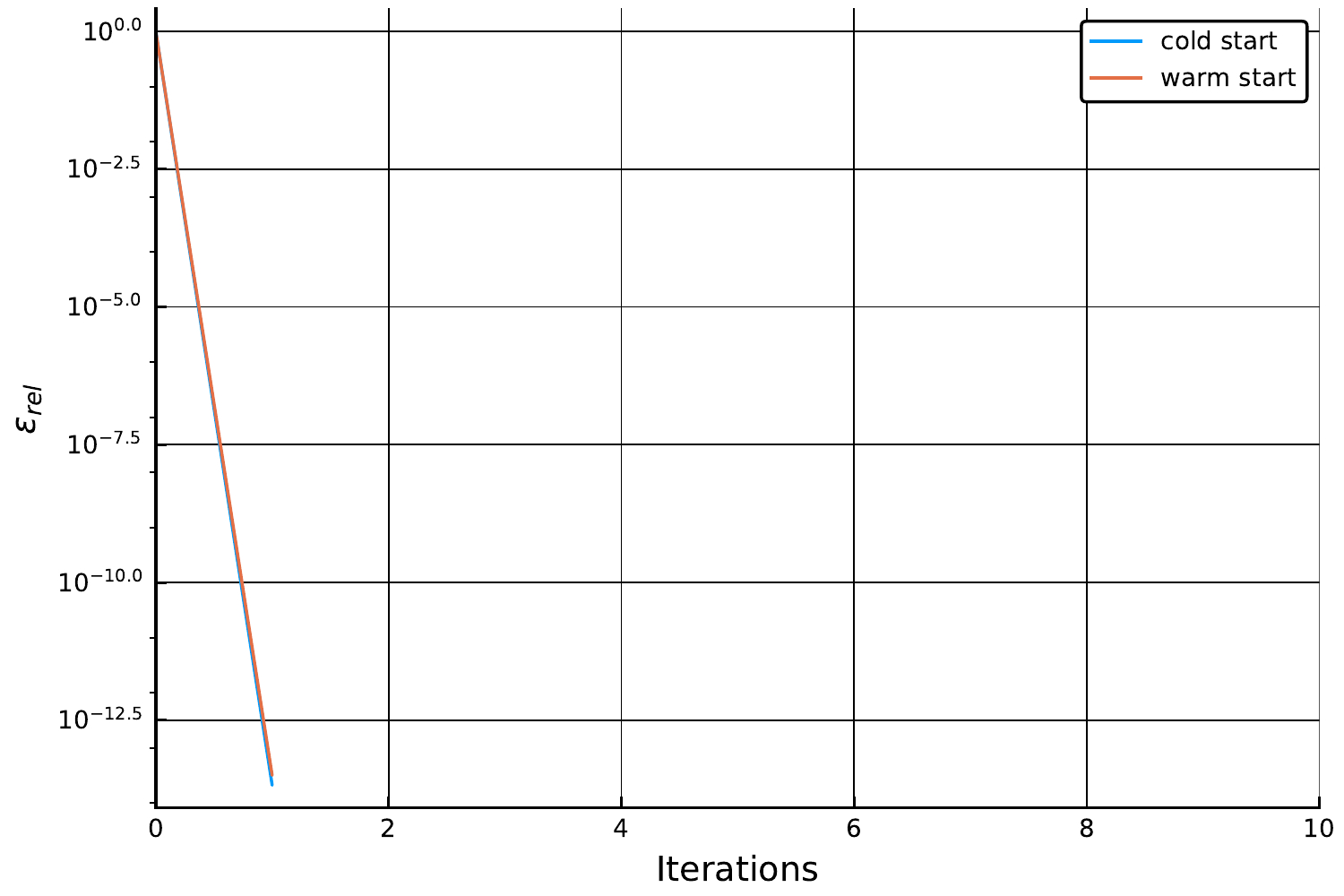"}
                \caption{Progression of the relative error}
        \end{subfigure}
        \hfill
        \begin{subfigure}{0.49\textwidth}
                \centering
                \includegraphics[width=\linewidth]{"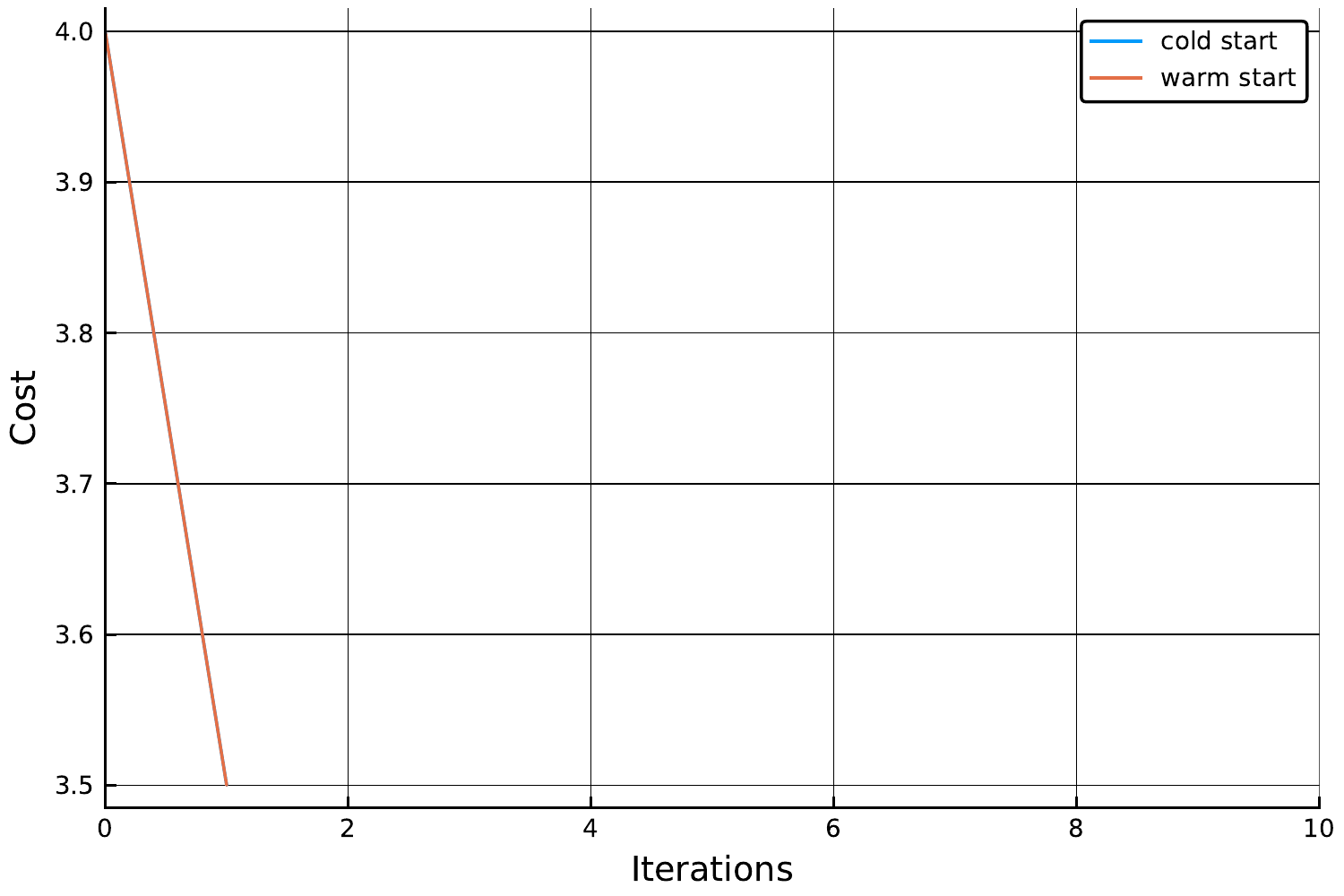"}
                \caption{Progression of the $\ell^2$-TV cost}
        \end{subfigure}\\
        \caption{Results of the proposed PD-RSSN method for the $\mathcal{P}(3)$-valued problem with known minimizer (compare \cref{fig: RSSN-known-minimizers-S2}). The method converges rapidly to the exact minimizer for both cold- and warm start.
        }
        \label{fig: RSSN-known-minimizers-SPD3}
 \end{figure}

\subsection{PD-RSSN for Solving Regularized $\ell^2$-TV} 
In this experiment we aim to compare the runtime performance of our PD-RSSN-based method to that of the   lRCPA algorithm at different accuracies. In contrast to the previous section, we focus on the more realistic 2D setting. We chose a dual regularization strength of $\beta = 10^{-6}$, as numerical experience showed that the generalized covariant derivative for both $S^2$ and $\mathcal{P}(3)$ can become singular without dual regularization.

For the $S^2$ problem, we used a $20\times 20$ artificial $S^2$ rotations image from \cite{bacak2016second} which is implemented in \manopt  with a half rotation 
around each axis. For the $\mathcal{P}(3)$ problem, we used a $10\times 10$ artificial $\mathcal{P}(3)$ image from \cite{bergmann2016parallel} which is also implemented in \manopt. 

Although lRCPA is not guaranteed to converge on positively-curved manifolds, it performed well in recent work \cite{bergmann2019fenchel}. The expectation is for lRCPA to be faster at the start but suffering from slow tail convergence. Hence, our method should give better performance for higher-accuracy solutions. 

\paragraph{Initialization}
For our numerical experiment, we measured the (CPU) runtime until the algorithms achieved $\epsilon_{rel} \in \{10^{-2},10^{-4},10^{-6}\}$. With a cold start, and even with a dual warm start, our method diverges, as can be expected from local methods. Therefore, we prepended lRCPA pre-steps with a coarse stopping criterion of $\epsilon_{rel} = 1/2$ and used the resulting primal-dual iterates to warm-start the RSSN solver. The computation of the relative errors includes the initial improvement due to the pre-steps with lRCPA.

\paragraph{Case 1: $\mathcal{M} = S^2$}
We computed an isotropically regularized ($q=2$) $\ell^2$-TV solution on the above data with $\alpha = 1.5$. We set $m_{i,j} := (0,0,1)^\top$ for all $ i,j=1,\ldots,20,$
so that again \jla{the tangent component of} $n$ is zero in $\mathcal{T}_m(S^2)^{20\times 20 \times 2}$. We initialized both PD-RSSN and lRCPA with $\sigma = \tau = 0.35$. For lRCPA we used $\gamma = 0.2$ to control the primal and dual step sizes.

\begin{figure}[h!]
        \centering
        \begin{subfigure}{0.49\textwidth}
                \centering
                 \includegraphics[width=0.9\linewidth]{"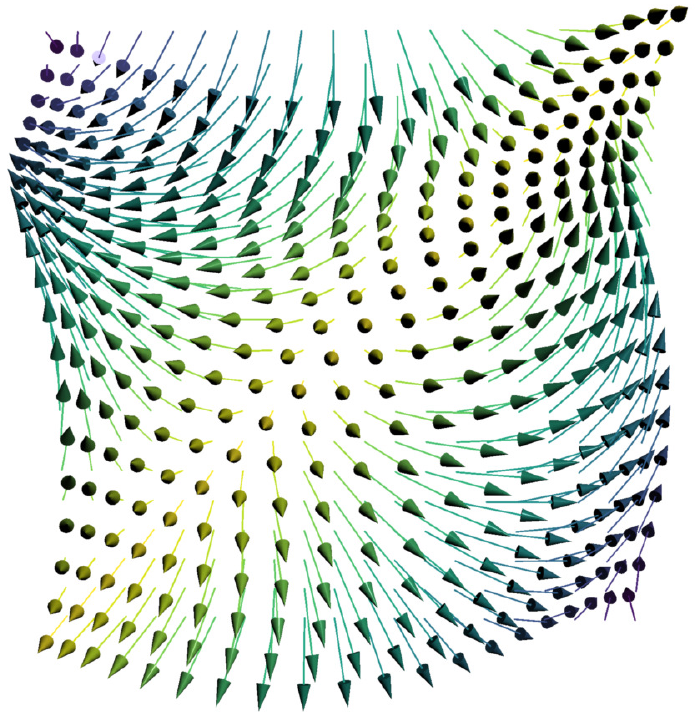"}
                \caption{Noisy input}
        \end{subfigure}
        \hfill
        \begin{subfigure}{0.49\textwidth}
                \centering
                 \includegraphics[width=0.9\linewidth]{"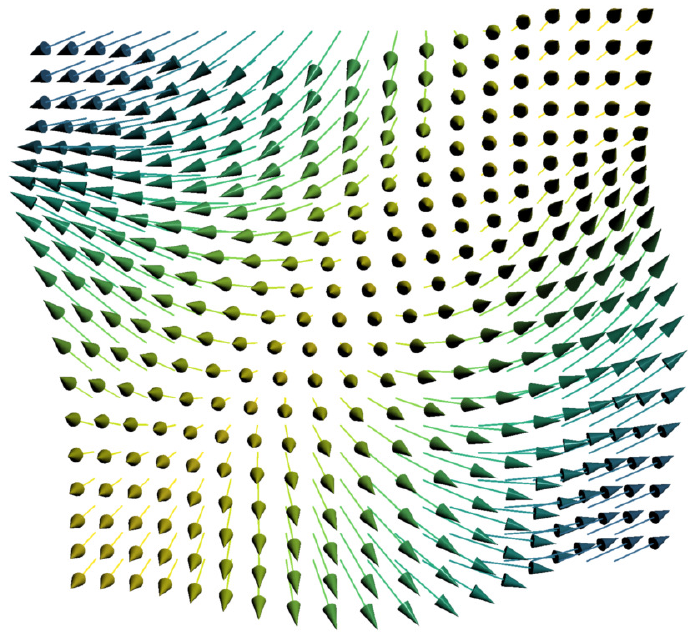"}
                \caption{Result}
        \end{subfigure}\\
        \begin{subfigure}{0.49\textwidth}
                \centering
                \includegraphics[width=\linewidth]{"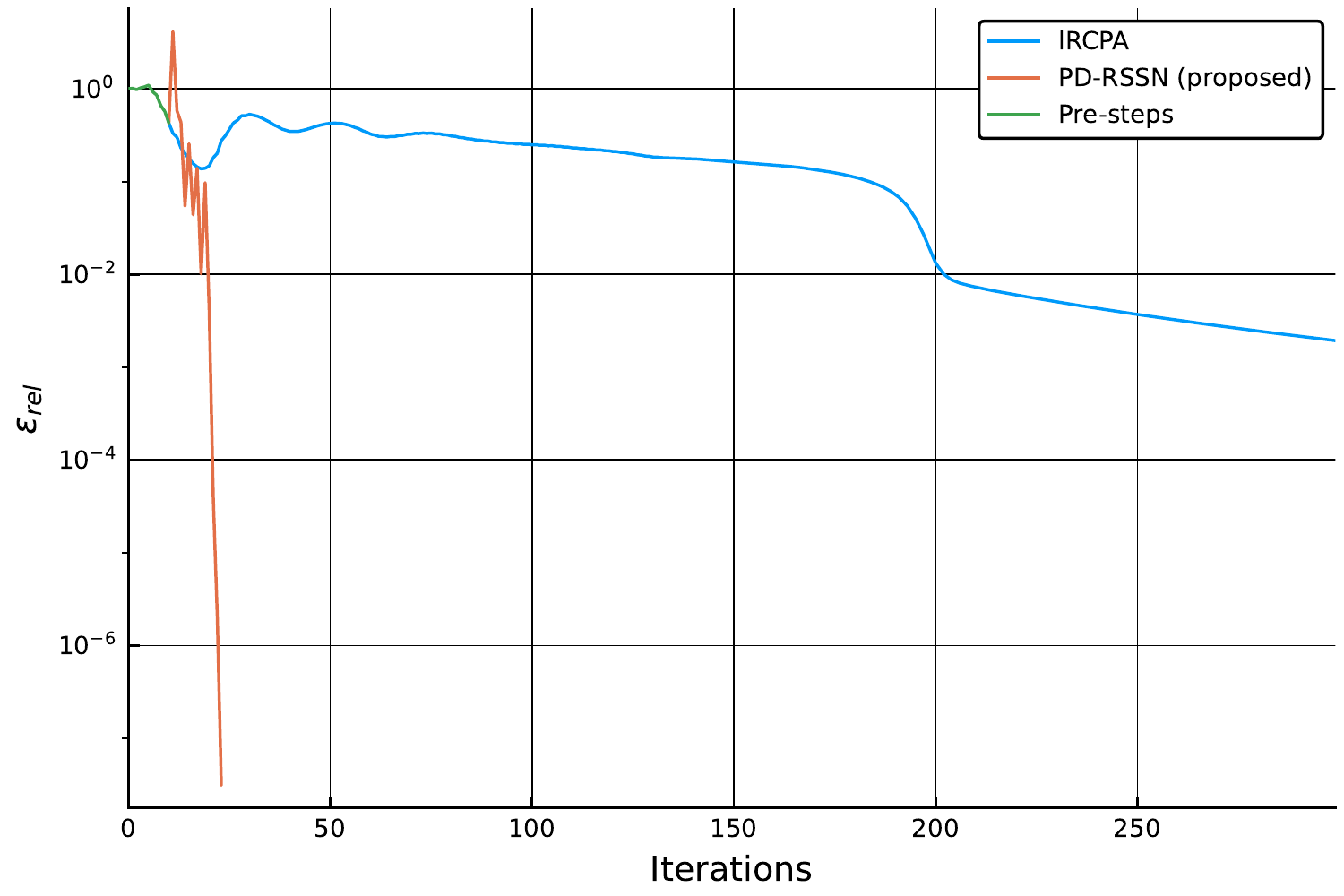"}
                \caption{Progression of the relative error}
                \label{fig:S2-relerror}
        \end{subfigure}
        \hfill
        \begin{subfigure}{0.49\textwidth}
                \centering
                \includegraphics[width=\linewidth]{"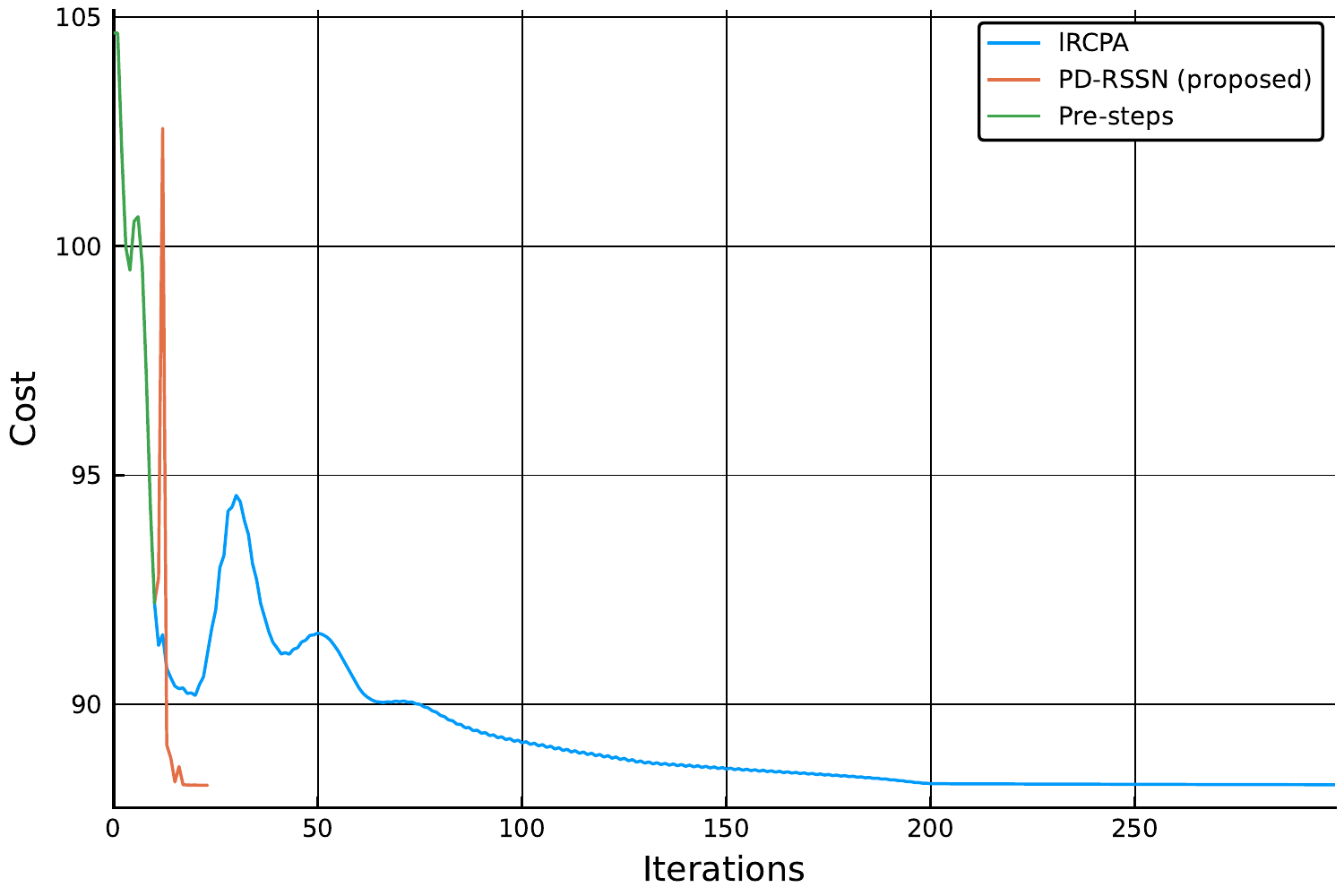"}
                \caption{Progression of the isotropic $\ell^2$-TV cost}
        \end{subfigure}
        \caption{Runtime behavior of lRCPA and the proposed higher-order PD-RSSN method  for an $S^2$-valued problem. The lRCPA method exhibits slow tail convergence as is typical for first-order methods. After few first-order pre-steps (shown in green), PD-RSSN  superlinearly converges within 13 iterations to an optimal solution.        }
        \label{fig: RSSN-comparison-of-algorithms 2}
 \end{figure}

The pre-steps took 2.437 seconds. The resulting runtimes are shown in \cref{tab: RSSN beta -6}. The solutions of lRCPA and PD-RSSN at $\epsilon_{rel} = 10^{-6}$ along with the progress of the relative error and the isotropic $\ell^2$-TV-cost are shown in \cref{fig: RSSN-comparison-of-algorithms 2}. 

\begin{table}[h!]
        \begin{tabular}{ccccccc}
                \toprule
                $S^2$  & \multicolumn{2}{c}{$\epsilon_{rel} = 10^{-2}$} & \multicolumn{2}{c}{$\epsilon_{rel} = 10^{-4}$} & \multicolumn{2}{c}{$\epsilon_{rel} = 10^{-6}$} \\ \cmidrule(l){2-3} \cmidrule(l){4-5} \cmidrule(l){6-7} 
                Method & Time (s)  & \# Iterations  & Time (s)  & \# Iterations   & Time (s)  & \# Iterations  \\ \midrule
                lRCPA  & \textbf{46.515} &  193  &  \textbf{159.11} &  697 & 608.781 & 2886  \\ 
                PD-RSSN  & 330.422 &  10 & 346.187 & 11  & \textbf{410.875} & 13  \\ \bottomrule
        \end{tabular}
        \caption{Runtime and iteration count for lRCPA (comparison) and the proposed PD-RSSN method to converge to different accuracies $\epsilon_{rel}$ for the two-dimensional $S^2$ problem. For high-accuracy solutions the proposed higher-order PD-RSSN method outperforms the first-order lRCPA algorithm with respect to runtime on a manifold with positive curvature.}
        \label{tab: RSSN beta -6}
\end{table}

The results confirm what we expect from a higher-order method: for higher accuracies with relative error $\epsilon_{rel}$ smaller than approximately $10^{-4}$, PD-RSSN cleary outperforms lRCPA. Furthermore, we see the superlinear convergence more clearly than in the 1D case due to the larger number of overall iterations\ (orange line in \cref{fig:S2-relerror}).

\paragraph{Case 2: $\mathcal{M} = \mathcal{P}(3)$}

Again we solved the isotropically regularized $\ell^2$-TV problem with $\alpha = 0.5$. We set choose $m_{i,j} := I$ for all $ i,j=1,\ldots,10$, 
so that $n$ is the zero vector \wdpa{in $\mathcal{T}_m\mathcal{P}(3)^{10\times 10 \times 2}$}. We initialized both PD-RSSN and lRCPA with $\sigma = \tau = 0.4$. For lRCPA we used $\gamma = 0.2$ to control the primal and dual step size.

\begin{figure}[h!]
        \centering
        \begin{subfigure}{0.49\textwidth}
                \centering
                \includegraphics[width=0.9\linewidth]{"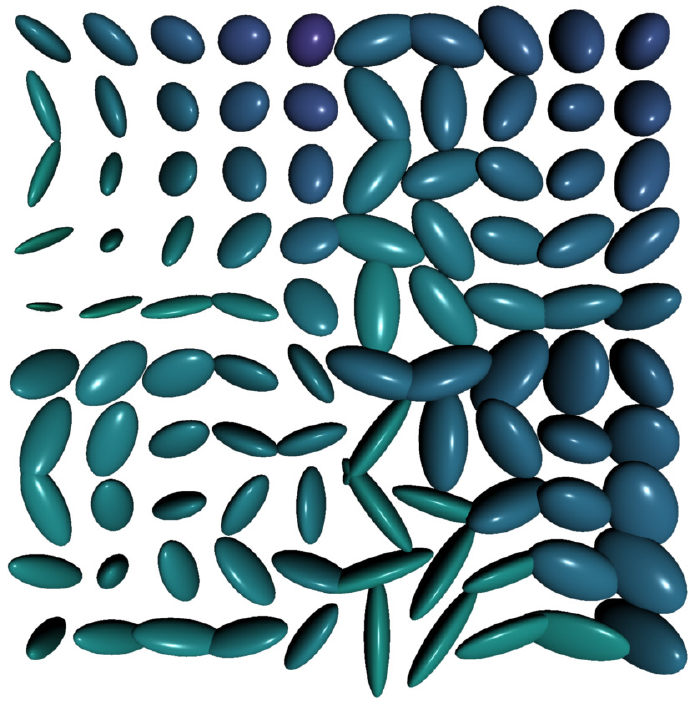"}
                \caption{Noisy input}
        \end{subfigure}
        \hfill
        \begin{subfigure}{0.49\textwidth}
                \centering
                \includegraphics[width=0.9\linewidth]{"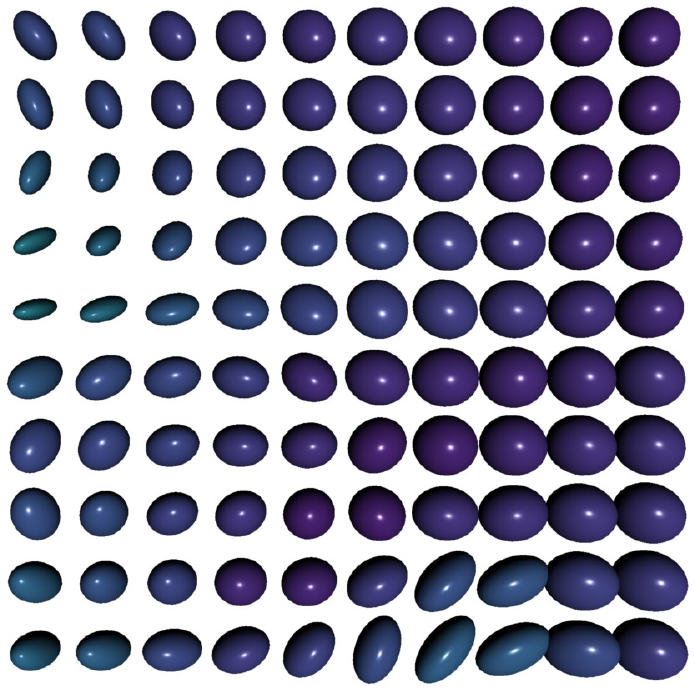"}
                \caption{Result}
        \end{subfigure}\\
        \begin{subfigure}{0.49\textwidth}
                \centering
                \includegraphics[width=\linewidth]{"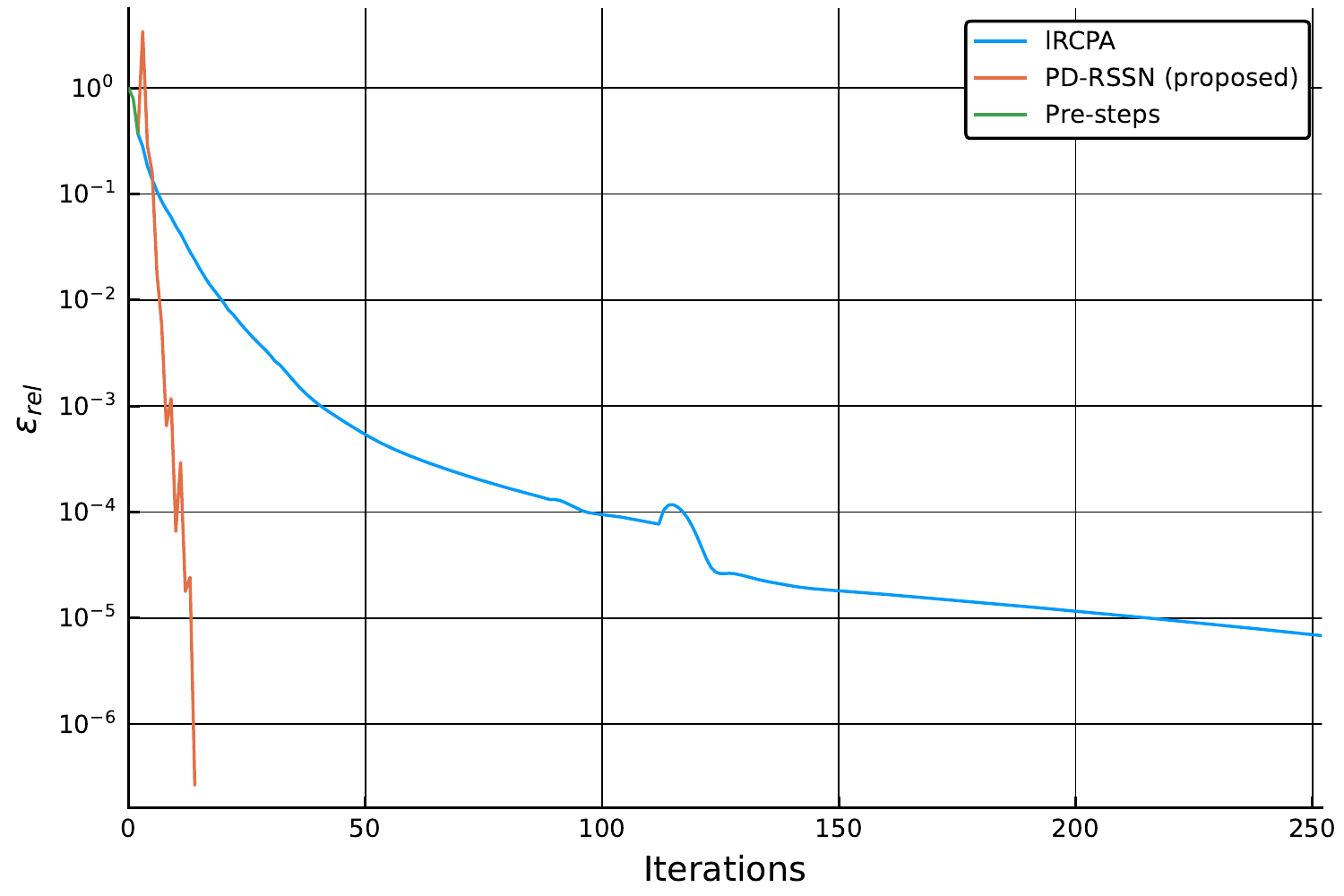"}
                \caption{Progression of the relative error}
        \end{subfigure}
        \hfill
        \begin{subfigure}{0.49\textwidth}
                \centering
                \includegraphics[width=\linewidth]{"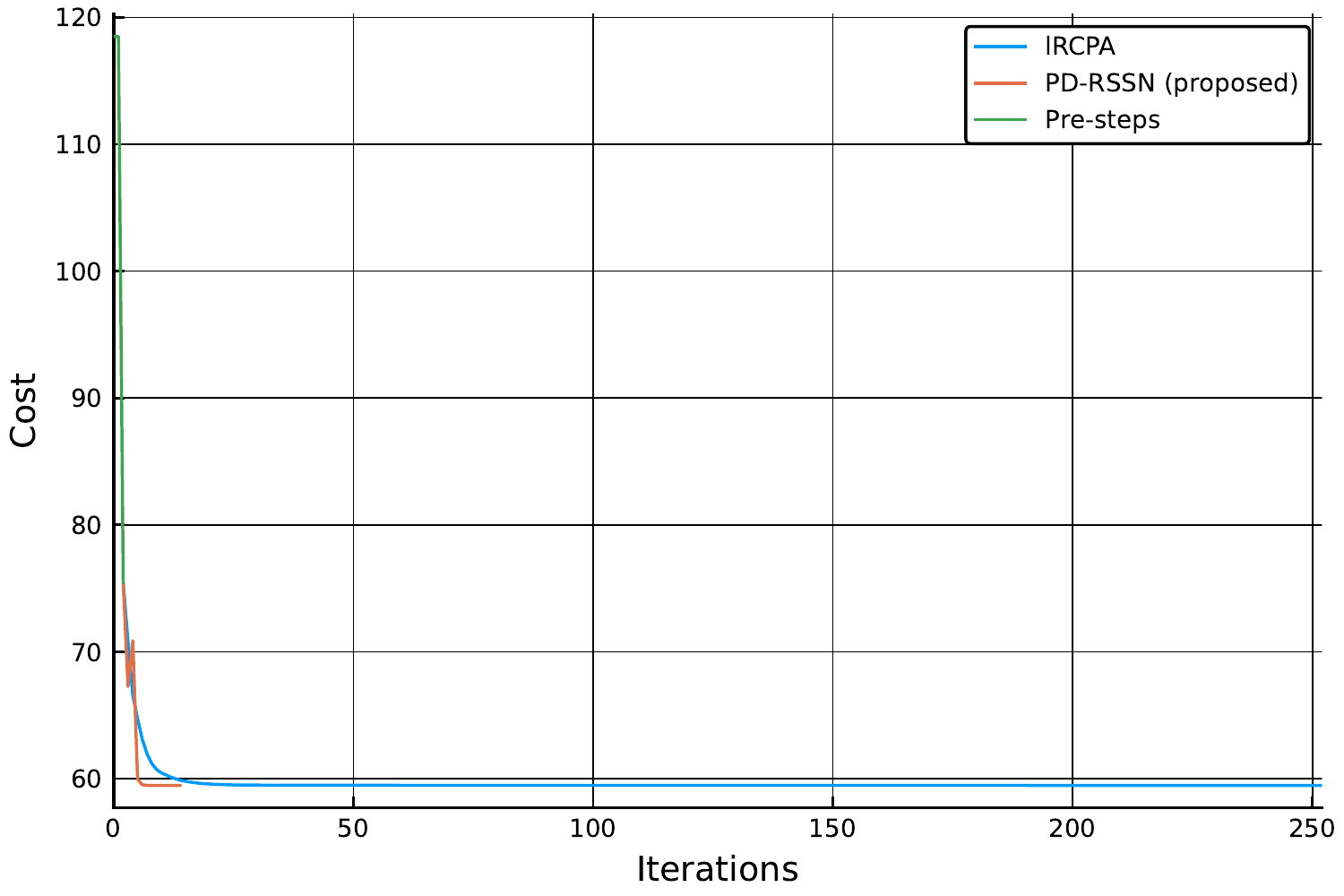"}
                \caption{Progression of the isotropic $\ell^2$-TV cost}
        \end{subfigure}
        \caption{Results on the negatively-curved manifold $\mathcal{P}(3)$ (compare \cref{fig: RSSN-comparison-of-algorithms 2});  first-order pre-steps are shown in green. The proposed PD-RSSN method superlinearly converges within 12 iterations after the pre-steps, while again lRCPA suffers from slow tail convergence.
        }
        \label{fig: RSSN-comparison-of-algorithms spd}
\end{figure}

The pre-steps took 9.718 seconds. The resulting runtimes are shown in \cref{tab: RSSN beta -6 SPD}. The solutions returned by lRCPA and PD-RSSN at $\epsilon_{rel} = 10^{-6}$ along with the error progress and the (isotropic) $\ell^2$-TV cost are shown in \cref{fig: RSSN-comparison-of-algorithms spd}. We note that lRCPA effectively stalled before reaching the relative error of $10^{-6}$ and was terminated after 2500 iterations.

\begin{table}[h!]
        \centering
        \begin{tabular}{ccccccc}
                \toprule
                $\mathcal{P}(3)$ &\multicolumn{2}{c}{$\epsilon_{rel} = 10^{-2}$} & \multicolumn{2}{c}{$\epsilon_{rel} = 10^{-4}$} & \multicolumn{2}{c}{$\epsilon_{rel} = 10^{-6}$} \\ \cmidrule(l){2-3} \cmidrule(l){4-5} \cmidrule(l){6-7} 
                Method &  Time (s) & \# Iterations  & Time (s)  & \# Iterations   & Time (s)  & \# Iterations  \\ \midrule
                lRCPA  &  \textbf{9.922} &  18  &  \textbf{48.36} &  97 & 1213.327 & $\geq$2500  \\
                PD-RSSN  &   237.062 &  5 & 380.047 & 8  & \textbf{556.328} & 12  \\ \bottomrule
        \end{tabular}
        \caption{Runtimes and iteration count for lRCPA (comparison) and the proposed PD-RSSN method to converge to different accuracies $\epsilon_{rel}$ for the two-dimensional  $\mathcal{P}(3)$ problem. As in the positively-curved $S^2$ case (\cref{tab: RSSN beta -6}), the proposed PD-RSSN algorithm outperforms lRCPA at higher accuracies on this negatively-curved manifold. lRCPA was terminated after 2500 iterations before achieving $\epsilon_{rel} = 10^{-6}$.}
        \label{tab: RSSN beta -6 SPD}
\end{table}

As for the $S^2$ example, for higher accuracies PD-RSSN clearly outperforms lRCPA. The latter even seems unable to reach such high accuracies in reasonable time, which aligns with the fact that first-order methods generally show sublinear convergence. For PD-RSSN, we again observe superlinear convergence.

\subsection{Convergence Behavior of PD-RSSN for Different Problem Sizes}
As computing a single Newton step scales quadratically at best, the convergence behavior of the proposed method for different problem sizes is important: Preferably, the number of required steps  should not, or only very weakly, dependent on the problem size. 

This type of convergence behavior is in general hard to measure precisely, because the problem will be slightly different at different scales. In order to get some intuition, we  performed the following experiment using data from the previous one: For $S^2$ we used the same $N\times N$ rotations image, and varied  $N = 10,15,20,25,30.$ For $\mathcal{P}(3)$, we used the same artificial $N\times N$ image with $N = 4,7,10,13,16$. We chose a tolerance of $10^{-6}$ throughout the experiment. The remaining parameters and the initialization were identical to the ones in the previous experiment.

The results are shown in \cref{tab: s2-scaling,tab: spd3-scaling}. We see that larger problem sizes typically require more iterations, but the dependence seems to be weak: While the problem size varied by approximately an order of magnitude, the number of required PD-RSSN steps remained in the range of 8 to 15.

We would like to note briefly that the number of pre-steps needed for the largest problem size in the $S^2$ example was considerably higher than for the other cases. Upon closer inspection, the convergence behavior in this case showed the typical oscillatory progression during the pre-steps as the lRCPA approach in \cref{fig:S2-relerror} (blue line). However, for the smaller problem sizes, the pre-step stopping criterion of $\epsilon_{rel}=1/2$ was achieved within the first downwards cycle. For the largest problem, this was not the case, which greatly increases the number of first-order steps until the pre-step stopping criterion is met (compare again the blue line in \cref{fig:S2-relerror}). 

\begin{table}[h!]
        \centering
        \begin{tabular}{ccccc}
                \toprule
                $S^2$ &\multicolumn{2}{c}{Pre-steps}  & \multicolumn{2}{c}{PD-RSSN} \\
                \cmidrule(l){2-3} \cmidrule(l){4-5}
                Problem size & Time (s)  & \# Iterations  & Time (s)  & \# Iterations     \\ \midrule
                100 & 0.237 & 5  &  6.043 &  8 \\
                225 & 0.519 & 7  &  26.417 &  10 \\
                400 & 0.990 & 10   &  138.523 &  13 \\
                625 & 2.147 & 14  &  577.480 &  15 \\
                900 & 42.321 & 220 &  1109.928 &  11 \\
                \bottomrule
        \end{tabular}
\caption{The effect of problem size on the amount of steps needed to converge to a solution with $\epsilon_{rel}=10^{-6}$ for the artificial problem on the $S^2$ manifold. The PD-RSSN steps seem to increase slightly for the first four problems, but drops for the largest problem. However, for this problem more pre-steps were needed to reach the tolerance $\epsilon_{rel} = 1/2$.}
\label{tab: s2-scaling}
\end{table}

\begin{table}[h!]
        \centering
        \begin{tabular}{ccccc}
                \toprule
                $\mathcal{P}(3)$ &\multicolumn{2}{c}{Pre-steps}  & \multicolumn{2}{c}{PD-RSSN} \\
                \cmidrule(l){2-3} \cmidrule(l){4-5}
                Problem size & Time (s) & \# Iterations  & Time (s)  & \# Iterations   \\ \midrule
                16 & 0.388 & 2 &  24.537 &  9 \\
                49 & 0.477 & 2  & 124.294 &  15 \\
                100 & 0.924 & 2 &  237.627 &  12 \\
                169 & 1.540 & 2 &703.031 &  15 \\
                256 & 3.465 & 3 & 1602.052 &  15 \\
                \bottomrule
        \end{tabular}
\caption{The effect of problem size on the amount of steps needed to converge to a solution with $\epsilon_{rel}=10^{-6}$ for the artificial problem on the $\mathcal{P}(3)$ manifold. Neither the pre-steps nor the PD-RSSN steps seem to be influenced much by the problem size.}
\label{tab: spd3-scaling}
\end{table}

\subsection{Primal-Dual Inexact Riemannian Semi-smooth Newton (PD-IRSSN)}
It is interesting to see if the convergence rates for the proposed inexact method in \cref{alg: inexact RSSN} as predicted by \cref{thm: inexact Newton} hold in practice. In particular, we expect linear convergence for constant relative residual  $a^k = \mathrm{constant}$ and superlinear convergence for $a^k \rightarrow 0$.  Given the similar behavior of the methods on positively and negatively curved manifolds in the previous experiments, we  only consider the $S^2$ manifold.

We used $\ell^2$-TV to denoise Bernoulli's Lemniscate \cite{bacak2016second}, which is a figure-8 on the 2-sphere. Note that this is once again a 1D problem. To this end, we sampled 128 $S^2$-valued points on the Lemniscate curve and distorted them with Gaussian noise with variance $\delta^2 = 0.01$. We employed the $\ell^2$-TV model with $\alpha = 0.5$ and $m$ set to the intersection point of the Lemniscate:
\begin{equation}
m_{i} := \frac{1}{\sqrt{2}}(1,0,1)^\top, \quad \text{ for } i=1,\ldots,128.
\end{equation}
We did not use dual regularization, i.e., 
$\beta = 0$.

\jla{Again we warm-started PD-IRSSN using lRCPA pre-steps} with $\sigma = \tau = 0.35,$ $\gamma = 0.2$, and  $\epsilon_{rel} = 10^{-1}$. In order to simulate inexact solution of the Newton steps, we defined three different progressions of the relative residual,
\begin{equation}
        a_1^k := 0, \qquad a_2^k := \frac{1}{5}, \qquad a_3^k := \frac{1}{5k}, 
\end{equation}
and randomly chose the corresponding absolute residual as
\begin{equation}
        r_i^k := a_i^k\|X(p^{k},\xi_n^k)\|_{(p^{k},\xi_n^k)}U_{(p^{k},\xi_n^k)}, \quad \text{ for } i = 1,2,3,
\end{equation}
where $U_{(p^{k},\xi_n^k)}$ is a tangent vector drawn from a normal distribution with unit covariance matrix. We then solved (to numerical precision) the modified inexact system for the step $d^k$:
\begin{equation}
        V_{(p^k,\xi_n^k)} d^k = X(p^k,\xi_n^k) + r_i^k.
\end{equation}
This allows to precisely control the residual and therefore simulate the behavior of PD-IRSSN at different accuracies.

\begin{figure}[h!]
        \centering
        \begin{subfigure}{0.49\textwidth}
                \centering
                \includegraphics[width=0.9\linewidth]{"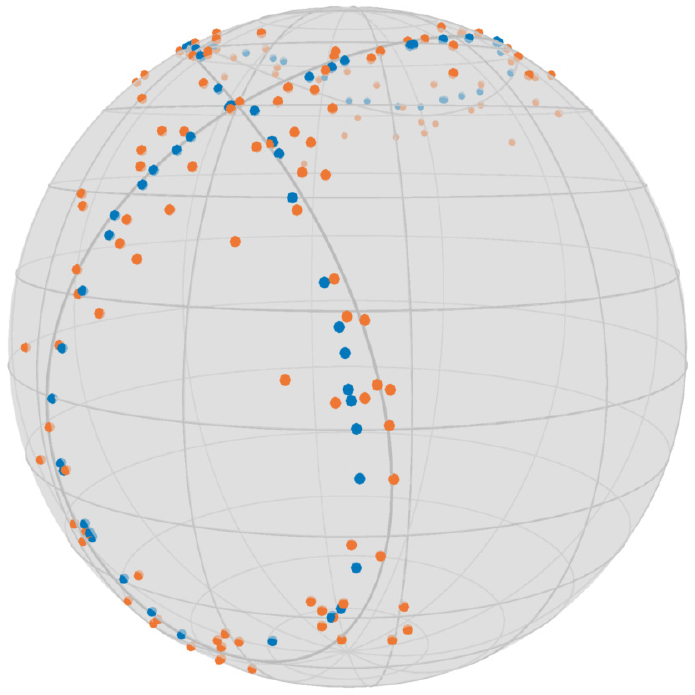"}
                \caption{Bernoulli's Lemniscate (grey), the noisy input (orange) and result (blue)}
        \end{subfigure}
 \hfill
 \begin{subfigure}{0.49\textwidth}
        \centering
        \includegraphics[width=\linewidth]{"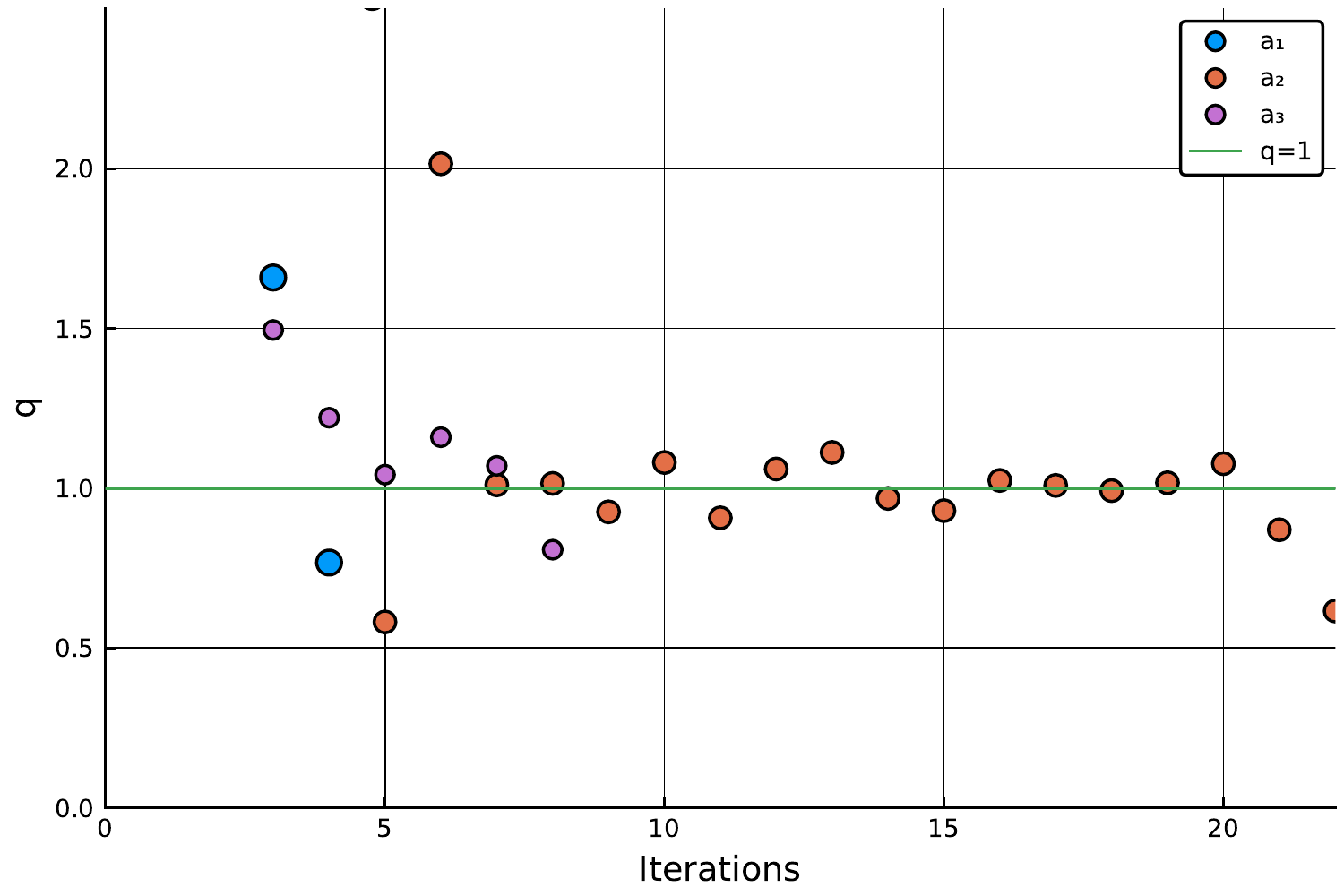"}
        \caption{Estimated convergence rates after the pre-steps}
 \end{subfigure}\\
        \begin{subfigure}{0.49\textwidth}
                \centering
                \includegraphics[width=\linewidth]{"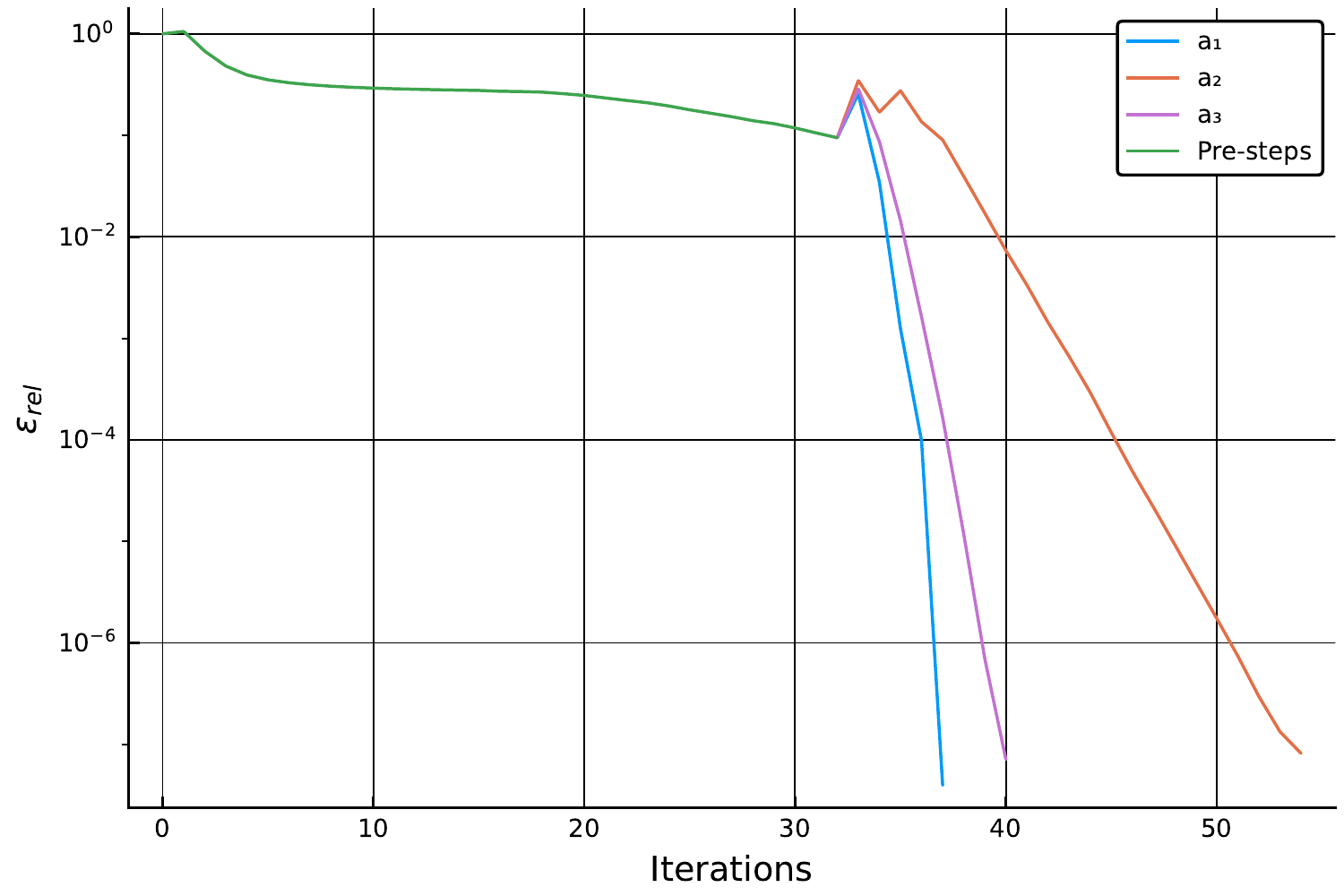"}
                \caption{Progression of the relative error of the lRCPA pre-steps (blue) and the three PD-IRSSN schemes}
        \end{subfigure}
        \hfill
        \begin{subfigure}{0.49\textwidth}
                \centering
                \includegraphics[width=\linewidth]{"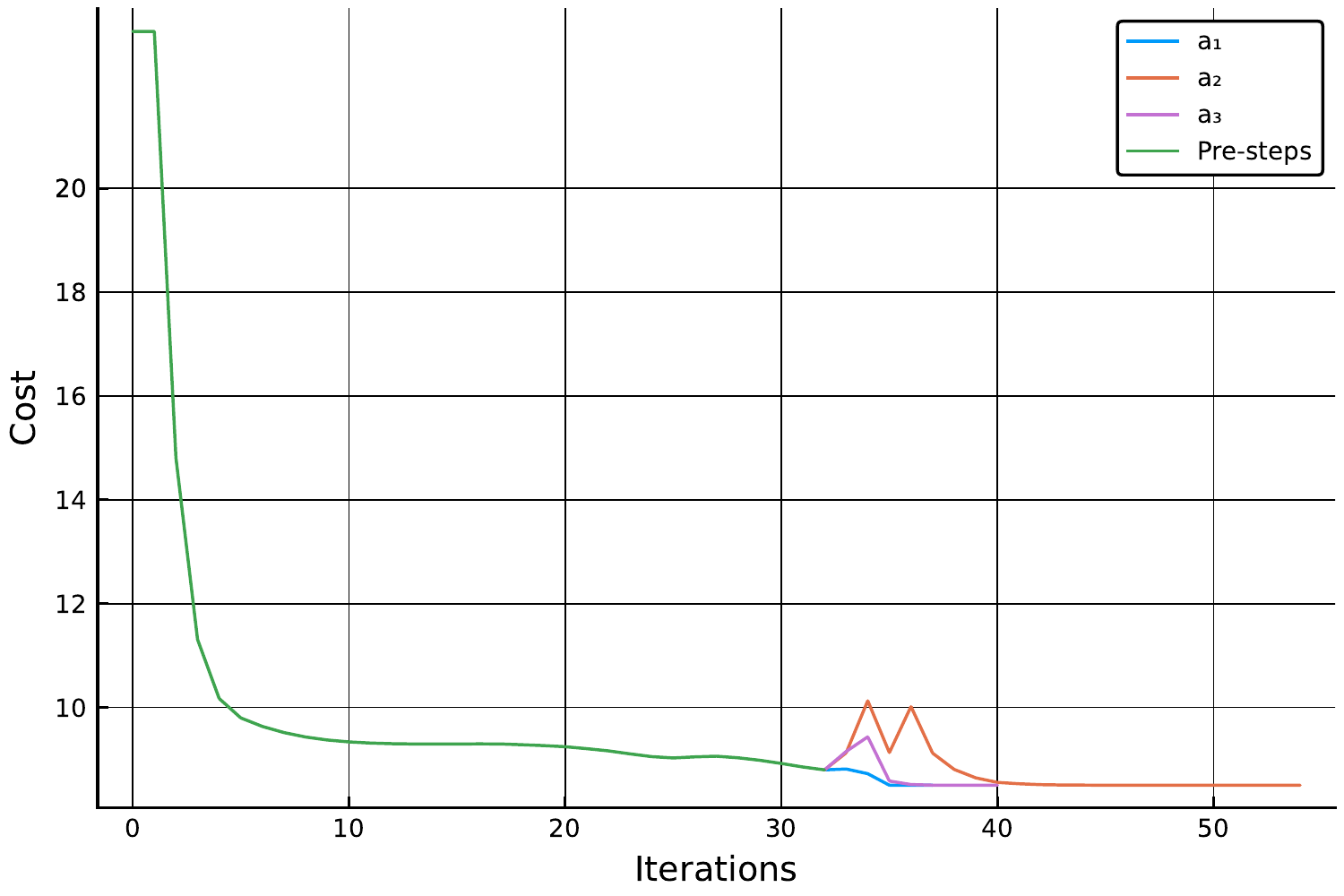"}
                \caption{Progression of the $\ell^2$-TV cost of the lRCPA pre-steps (blue) and the three PD-IRSSN schemes}
        \end{subfigure}
        \caption{Convergence behavior of the Primal-Dual Inexact Riemannian Semi-smooth Newton method (PD-IRSSN) for recovering Bernoulli's Lemniscate from noisy input data. As predicted by \cref{thm: inexact Newton},  exact PD-IRSSN, i.e., PD-RSSN, ($a_1$) and PD-IRSSN with progressively increasing accuracy ($a_3$) converge superlinearly. For constant relative accuracy ($a_2$), we still observe at least linear convergence. 
        }
        \label{fig: RSSN-inexact-semismooth}
 \end{figure}

From the error progression, we estimated a convergence rate $q$ by
\begin{equation}
        q^k:= \frac{\log \left(\frac{\|X(p^{k},\xi_n^{k})\|}{\|X(p^{k-1},\xi_n^{k-1})\|}\right)}{\log \left(\frac{\|X(p^{k-1},\xi_n^{k-1})\|}{\|X(p^{k-2},\xi_n^{k-2})\|}\right)}.
\end{equation}
The results are shown in \cref{fig: RSSN-inexact-semismooth}. As predicted by \cref{thm: inexact Newton}, for the exact method with zero residual $a_1$ we obtain clear superlinear convergence. For constant relative residual $a_2$, convergence is clearly linear with rate $q =1$, and  for decreasing relative residual $a_3$,  convergence qualitatively appears superlinear -- although hard to judge due to the rapid convergence -- in agreement with the theorem.

\section{Conclusions}
\label{sec:conclusions}


In this work we have proposed the Primal-Dual Riemannian Semi-smooth Newton (PD-RSSN) method as a higher-order alternative for solving non-smooth variational problems on manifolds. Although \wdpa{the method is only} locally convergent, it experimentally performs well on  isotropic and anisotropic 1D and 2D $\ell^2$-TV problems. On both $S^2$ and $\mathcal{P}(3)$, representing manifolds with positive and negative curvature, we observe superlinear convergence, which allows to generate high-accuracy solutions. In particular, PD-RSSN can be used to overcome the slow tail convergence that is a common bottleneck of first-order methods. 

Moreover, we have proposed an inexact variant of RSSN (IRSSN) with provable local convergence rates. The theoretical results are supported by numerical experiments, in which we observe at least linear convergence in the inexact case, and superlinear convergence when accuracy is progressively increasing.

%
%
%

\section*{Acknowledgments}
We would like to thank Ronny Bergmann for fruitful discussions regarding differential-geometric interpretations of earlier work and his Julia library \manopt. 

\clearpage

\bibliographystyle{siamplain}
\bibliography{references}

\end{document}